\documentclass{amsart}
\usepackage{amsmath, amsthm, diagrams, latexsym, amssymb, amsfonts, epsfig, epsf}
\usepackage[dvips]{color}

\newenvironment{pf}{\proof[\proofname]}{\endproof}
\theoremstyle{plain}
\newtheorem{Th}{Theorem}[section]
\newtheorem{Cor}[Th]{Corollary}
\newtheorem{Prop}[Th]{Proposition}
\newtheorem{Lemma}[Th]{Lemma}
\numberwithin{equation}{section} 
\numberwithin{figure}{section} 
\theoremstyle{definition}
\newtheorem{Rem}[Th]{Remark}

\newtheorem{Example}[Th]{Example}
\newtheorem{Def}[Th]{Definition}

\newarrow{Dashto}{}{dash}{}{dash}{>}

\newcommand{\cal}[1]{\mathcal{#1}}
\newcommand{\C}{\mathbb C}

\newcommand{\Z}{\mathbb Z}
\newcommand{\R}{\mathbb R}
\newcommand{\G}{\Gamma}

\newcommand{\D}{\Delta}

\newcommand{\cB}{\cal B}

\newcommand{\cF}{\cal F}

\newcommand{\cJ}{\cal J}
\newcommand{\cM}{\cal M}
\newcommand{\cO}{\cal O}
\newcommand{\cU}{\cal U}

\newcommand{\p}{\partial}
\newcommand{\T}{{\mathbb T}}
\newcommand{\e}{\varepsilon}

\newcommand{\sig}{\sigma}
\newcommand{\Sig}{\Sigma}

\newcommand{\Res}{\operatorname{Res}}

\newcommand{\sgn}{\operatorname{sgn}}

\newcommand{\inter}{\operatorname{int}}
\newcommand{\Vertex}{\operatorname{Vert}}
\newcommand{\Mat}{\operatorname{Mat}}

\newcommand{\cdeg}{\operatorname{cdeg}}

\newcommand{\Tr}{\operatorname{Tr}}

\newcommand{\rs}[1]{Section~\ref{S:#1}}
\newcommand{\rl}[1]{Lemma~\ref{L:#1}}
\newcommand{\rp}[1]{Proposition~\ref{P:#1}}

\newcommand{\rex}[1]{Example~\ref{Ex:#1}}
\newcommand{\re}[1]{(\ref{e:#1})}
\newcommand{\rc}[1]{Corollary~\ref{C:#1}}
\newcommand{\rt}[1] {Theorem~\ref{T:#1}}
\newcommand{\rd}[1]{Definition~\ref{D:#1}}
\newcommand{\rf}[1]{Figure~\ref{F:#1}}

\begin{document}


\title{Combinatorial construction of toric residues}
\author[Amit Khetan]{Amit Khetan}
   \address[Amit Khetan]{Department of Mathematics and Statistics\\
University of Massachusetts\\ Amherst, MA 01003, USA}
   \thanks{Amit Khetan was supported by NSF postdoctoral fellowship DMS-0303292.}
   \email{khetan@math.umass.edu}
   \author[Ivan Soprounov]{Ivan Soprounov}
   \address[Ivan Soprounov]{Department of Mathematics and Statistics\\
University of Massachusetts\\ Amherst, MA 01003, USA}
   \email{isoprou@math.umass.edu}
\keywords{Toric varieties, toric residues, semi-ample degrees, facet colorings,
combinatorial degree}
\subjclass[2000]{Primary 14M25; Secondary 52B20, 06A07}
   

\begin{abstract} 
The toric residue is a map depending on $n+1$ divisors on a
complete toric variety of dimension $n$. It appears in a variety of
contexts such as sparse polynomial systems, mirror symmetry, and GKZ
hypergeometric functions.  In this paper we investigate the problem of
finding an explicit element whose toric residue is equal to one. Such an
element is shown to exist if and only if the associated polytopes are
essential. We reduce the problem to finding a collection of partitions
of the lattice points in the polytopes satisfying a certain
combinatorial property. We use this description to solve the problem
when $n=2$ and for any $n$ when the polytopes of the divisors share a
complete flag of faces. The latter generalizes earlier results when
the divisors were all ample.
\end{abstract}

\maketitle

\section{Introduction}

Toric residues are fundamental invariants of sparse polynomial
systems. They were first studied by Cox \cite{Coxres} who defined the
residue of $n+1$ sections of an ample line bundle on a toric variety
$X$. The definition was extended by Cattani, Cox, and Dickenstein to
sections of $n+1$ arbitrary line bundles \cite{CCD}. There are numerous
applications to sparse resultants and resultant or subresultant
complexes \cite{CDS, AK}, mixed Hodge structures \cite{BC}, and mirror symmetry
\cite{BM}.

The related notion of global residue in the torus, a sum of
Grothendiek local residues, was studied by Gelfond, Khovanskii, and
Soprounov \cite{GKh2, Sop1}. Cattani, Cox and Dickenstein \cite{CCD} showed that the global
residue could always be computed as an instance of the toric residue.
Applications of the toric and global residue include GKZ hypergeometric systems \cite{CaD2, CDS3, CDS4}
and computations on sparse polynomial systems such as counting the
number of real roots and computing elementary symmetric functions on
the roots \cite{CDS2, GKh2}.

Given $n+1$ arbitrary sparse Laurent polynomials $f_0, \dots, f_n$ in
$n$ affine variables, let $P_0, \dots, P_n$ be their corresponding
Newton polytopes.  The Minkowski sum $P = P_0 +\dots+ P_n$ determines a
toric variety $X$, and each $P_i$ corresponds to a semi-ample divisor
class. In the homogeneous coordinate ring $S$ of $X$ each $f_i$ can be
homogenized to a polynomial $F_i$ of degree $\alpha_i$ corresponding
to the divisor class of $P_i$. The toric residue $\Res_F$ is a linear
function on homogeneous polynomials of a certain critical degree corresponding
to the interior of $P$ which vanishes on the ideal of the $F_i$.

In many cases of interest, for example when all of the $P_i$ are full
dimensional, the ideal of the $F_i$ has codimension $1$ in the
critical degree. Hence knowing a single element of non-zero residue
will allow a full computation of the residue map. More generally, we
show in \rs{essential} that there is an element of non-zero
residue whenever the polytopes form an essential family. 
The goal of this paper
is a general framework for the construction of specific elements
whose residue we can compute. The construction depends only on the
combinatorics and affine geometry of the polytopes $P_i$.

\begin{Th}\label{T:main1} Let $X$ be a complete toric variety of dimension $n$. 
Fix $n+1$ semi-ample degrees $\alpha_0,\dots,\alpha_n$ on $X$ and let
$P_0,\dots,P_n$ be their polytopes.  Let
$$P_i\cap\Z^n=M_{i0}\sqcup\dots\sqcup M_{in},\quad 0\leq i\leq n,$$
be a  collection of partitions of the lattice points of the $P_i$
such that
\begin{enumerate}
\item for any lattice point $u\in M_{ij}$, at least one vertex of the minimal
face of $P_i$ containing $u$ lies in $M_{ij}$,
\item for any permutation $\e$ of $\{0, \dots, n\}$:
$$\sum_{i=0}^{n}M_{\e(i)i} \subset \inter\Big(\sum_{i=0}^nP_i\Big).$$
\end{enumerate}
Given a collection of Laurent polynomials $f_0,\dots,f_n$ supported on $P_0,\dots,P_n$
$$f_i=\sum_{u\in P_i\cap\Z^n}c_ut^u,\quad 0\leq i\leq n$$
define polynomials
$$f_{ij}=\sum_{u\in M_{ij}}c_ut^u,\quad 0\leq i,j\leq n.$$ Then
$h=\det(f_{ij})$ is a Laurent polynomial supported on
$\inter(\sum_{i=0}^nP_i)$. The toric residue $\Res_F(H)$ of the corresponding homogeneous polynomial
$H$ of critical degree for the homogenized $F_0,\dots,F_n$
is an integer that depends only on the
combinatorics of the $P_i$ and the partitions of their lattice points.
\end{Th}

Using this theorem we are able to find an element of residue $\pm
1$, i.e. find an appropriate collection of partitions, in two
important cases.  The first is when $P_i$ share a complete flag of
faces. This will generalize earlier results of D'Andrea and Khetan
when all of the $\alpha_i$ were ample degrees.  The second application
is a complete analysis when $n=2$. We show that, except for one
degenerate family of supports, we can always find a collection of
partitions yielding an element of residue $\pm 1$.

The proof of the theorem makes use of  some very elegant
combinatorics. Starting with a partition of the lattice points we
will show that there are induced colorings of the faces of the polytope
$P = \sum P_i$. Moreover, the matrix will yield a canonical coloring
of the facets of the barycentric refinement of $P$. Such a facet
coloring will allow us to reduce the computation to that of the
residue of a monomial with respect to a monomial ideal. By an
earlier theorem of Soprounov \cite{Sop}, the residue is the {\em combinatorial
degree} of the coloring which can be computed by counting the number
of flags of certain colors.

The paper is organized as follows.  \rs{prelim} provides the
definitions of the toric residue and some basic
properties. \rs{essential} proves the existence of elements of
non-zero residue if and only if the polytopes are
essential. \rs{cdeg} introduces facet colorings of polytopes and their
connection to the toric residue of monomials. The residue for general polynomials
is reduced to the monomial case via the Global Transformation Law. 
\rs{partition} and \rs{canoncolor}
discuss the relationships between partitions, colorings, and residue
matrices used to complete the proof of \rt{main1}. 
\rs{locallyunmixed} uses the previous results to give an explicit
element of residue 1 when the polytopes $P_i$ share a complete flag of
faces.  \rs{dim2} is a complete analysis when $X$ is of dimension 2.
Finally, \rs{conclusion} discusses progress in dimensions three and
higher.

\section{Preliminaries} \label{S:prelim}

We begin by setting up the notation and reviewing some basic definitions and facts about
toric varieties and toric residues. For details and proofs we refer the reader to \cite{CCD, Coxhom, Coxres, F}.

\subsection{Toric residue} Consider an $n$-dimensional complete toric variety $X$
determined by a rational complete fan $\Sig\subset\R^n$. 
Let $\Sig(1)$ denote the set of 1-dimensional cones (rays) of $\Sig$.
Each ray $\rho\in\Sig(1)$ determines a $\T$-invariant
irreducible divisor $D_\rho$ on $X$. As introduced by Cox in \cite{Coxhom} the variety $X$ has the
homogeneous coordinate ring $S=\C[x_\rho : \rho\in\Sig(1)]$ graded by the
Chow group $A_{n-1}(X)$ so that a monomial
$x^a=\prod_{\rho}x_\rho^{a_\rho}$ has degree
$$\deg(x^a)=\big[\sum_{\rho\in\Sig(1)}a_{\rho}D_{\rho}\big]\in A_{n-1}(X).$$
Denote by $S_\alpha$ the graded piece of $S$ consisting of all
polynomials of degree $\alpha\in A_{n-1}(X)$. 

Let $D=\sum_{\rho}a_\rho D_\rho$ be a representative of $\alpha\in A_{n-1}(X)$.
It defines a continuous piecewise linear function $\psi_D$ on the support $|\Sig|$ such that
$\psi_D(v_\rho)=-a_\rho$ for all $\rho\in\Sig(1)$, where $v_\rho$ denotes the primitive generator
of $\rho$ (see \cite[Section 3.3]{F}). It also determines a convex polytope
$$P_D=\{u\in\R^n\,:\,\langle u,v_\rho\rangle\geq-a_\rho,\ \rho\in\Sig(1)\}=\{u\in\R^n\,:\,u\geq \psi_{D} \text{ on } |\Sig|\}.$$
To every lattice point $u$ of $P_D$ we can assign a monomial $\chi^u$ in $S$ of degree $\alpha$:
$$\chi^u=\prod_{\rho\in\Sig(1)}x_\rho^{\langle u,v_\rho\rangle+a_\rho},\quad u\in P_D\cap\Z^n.$$
One can check that this map is a bijection. Furthermore, 
given a Laurent polynomial $f(t)=\sum_u c_u t^u$ supported in $P_D$ its {\it $P_D$-homogenization} 
is the homogeneous polynomial
\begin{equation}\label{e:homocoord}
F=\sum_{u\in P_D\cap\Z^n}c_u\chi^u=\sum_{u\in P_D\cap\Z^n}c_u\prod_{\rho\in\Sig(1)}x_\rho^{\langle u,v_\rho\rangle+a_\rho}\in S_\alpha.
\end{equation}
Notice that if $f$ is supported in the interior of $P_D$ then the $P_D$-homogenization is divisible
by the product of all the variables $x_\rho$, $\rho\in\Sig(1)$.
It is easy to see that if $D$ and $D'$ are linearly equivalent then $\psi_D-\psi_{D'}$ is a linear function,
and $P_D$ and $P_{D'}$ are the same up to a translation. Therefore, $P_D$-homogenization is independent
of the choice of the representative $D$ of the divisor class $\alpha$. In what follows the {\it polytope of $\alpha$}
will mean the polytope of any representative of $\alpha$ and will be denoted by $P_\alpha$.

Recall the construction of the Euler form $\Omega$ from \cite{CCD}. Let $(e_1,\dots,e_n)$ be
a basis for $\Z^n$ and for every subset $I\subset\Sig(1)$ of size $n$ denote
$$\det(\eta_I)=\det(\langle e_i,v_{\rho}\rangle\,:\,1\leq i\leq n,\,\rho\in I),
\quad dx_I=\wedge_{\rho\in I}\,dx_{\rho},\quad \hat x_{I}=\prod_{\rho\not\in I}x_\rho.$$
Then the {\it Euler form} on $X$ is the sum over all size $n$ subsets $I\subset\Sig(1)$:
$$\Omega=\sum_{|I|=n}\det(\eta_I)\hat x_{I}dx_I.$$

Now we recall the definition of the toric residue \cite{Coxres, CCD}.
Consider $n+1$ homogeneous polynomials $F_i\in S_{\alpha_i}$, for $0\leq i\leq n$.
Their critical degree is defined to be
$$\nu=\sum_{i=0}^n\alpha_i-\sum_{\rho}\deg(x_\rho).$$
Then for every polynomial $H$ of degree $\nu$ consider a meromorphic $n$-form on $X$:
$$\omega_F(H)=\frac{H\Omega}{F_0\cdots F_n},$$
where $\Omega$ is the Euler form. We use $F$ to denote the list $(F_0,\dots,F_n)$.
Suppose that the~$F_i$ do not vanish simultaneously on $X$. Then $X$
has an open cover $\cU$ by the $n+1$ sets $U_i=\{x\in X : F_i(x)\neq 0\}$ 
and $\omega_F(H)$ defines a \v Cech cohomology class $[\omega_F(H)]\in H^n(X,\widehat\Omega_X^n)$
relative to the cover $\cU$. Here $\widehat\Omega_X^n$ denotes the sheaf of Zariski $n$-forms on $X$. 
One can check that the class $[\omega_F(H)]$ is
alternating in the order of the $F_i$ and is zero if $H$ belongs
to the ideal of $F_0,\dots,F_n$. Therefore, $[\omega_F(H)]$ depends on the equivalence
class of $H$ modulo the ideal $\langle F_0,\dots,F_n\rangle$. The {\it toric
residue map}
$$\Res^X_F:S_\nu/\langle F_0,\dots,F_n\rangle_\nu\to\C,$$
is given by $\Res^X_F(H)=\Tr_X([\omega_F(H)])$, where $\Tr_X$ is the
trace map on $X$. When there is no danger of confusion we will write $\Res_F(H)$ instead of $\Res^X_F(H)$.

\subsection{Semi-ample degrees} Let $X$ be a complete $n$-dimensional toric variety defined by a complete
fan $\Sig$ in $\R^n$. Recall that a $\T$-Cartier divisor $D$ on $X$ is called
{\it semi-ample} if the corresponding line bundle $\cO(D)$ is generated by global sections.
Equivalently, $D$ is semi-ample if and only if the corresponding piecewise linear function
$\psi_D$ is convex \cite[Section 3.4]{F}. Consider the (generalized) normal fan $\Sig_D$ of the
polytope $P_D$ of $D$, i.e. a complete fan whose cones are
$$\sig_\G=\{v\in(\R^n)^* : \langle u,v\rangle\geq \langle u',v\rangle,\text{ for all } u\in P_D,u'\in \G\},$$
for every face $\G$ of $P_D$. It follows that if $D$ is semi-ample then
$\Sig$ refines $\Sig_D$. Indeed, by the convexity of $\psi_D$
for any maximal cone $\sig\in\Sig$ the restriction of $\psi_D$ to $\sig$ defines a
a vertex $u$ of $P_D$. Then $\sig\subset\sig_u$, $\sig_u\in\Sig_D$.
We will say that a degree $\alpha=[D]\in A_{n-1}(X)$ is {\it semi-ample} if $D$ is semi-ample.

Consider a collection of $n+1$ semi-ample degrees $\alpha_0,\dots,\alpha_n$ on $X$. Let $P_0,\dots,P_n$
be their polytopes (defined up to translations) and $\Sig_0,\dots,\Sig_n$ the normal fans of the polytopes.
By above $\Sig$ refines each $\Sig_i$ and, thus, refines the minimal common refinement of the $\Sig_i$, which is
the normal fan $\Sig_P$ of the Minkowski sum $P=\sum_{i=0}^nP_i$ by \cite[Chapter 5, Theorem 4.8]{GKZ}.

Now let $\pi: X'\to X$ be a birational morphism  defined by a refinement $\Sig'\to\Sig$. 
If $D$ is a $\T$-Cartier divisor on $X$ then the pull-back $\pi^*(D)$ has the same
piecewise linear function $\psi_D$ and the same polytope $P_D$. It follows that 
$\alpha'=\pi^*(\alpha)$ is semi-ample on $ X'$ if $\alpha$ is semi-ample on $X$. 
Also if $F$ is a homogeneous polynomial in $S_\alpha$ and $f$ the corresponding
Laurent polynomial supported in $P_\alpha$ then the pull-back $ F'=\pi^*(F)$ is
the $P_\alpha$-homogenization of $f$ in the homogeneous coordinate ring $S'$ of $ X'$, and
hence $ F'\in S'_{\alpha'}$.

Next we will see how the toric residue $\Res^X_F$ behaves under the birational morphism $\pi: X'\to X$.

\begin{Prop}\label{P:birational} Let $X$ be a complete $n$-dimensional toric variety defined by a complete fan $\Sig$.
Let $\pi: X'\to X$ be a birational morphism induced by a refinement $\Sig'\to\Sig$.
Suppose $\alpha_0,\dots,\alpha_n$ are semi-ample degrees with polytopes $P_0,\dots,P_n$ and consider $n+1$ polynomials
$F_i\in S_{\alpha_i}$ not vanishing simultaneously on $X$. Then the polynomials
$ F'_i=\pi^*(F_i)\in S'_{\alpha'_i}$ do not vanish simultaneously on $ X'$.
Furthermore, let $g$ be any Laurent polynomial supported in the interior of $P=\sum_{i=0}^n P_i$, and
$G$ (resp. $ G'$) be the $P$-homogenization of $g$ in $S$ (resp. $ S'$). Then
the homogeneous polynomials $H=G/x_{\Sig(1)}$ and 
$ H'= G'/x_{\Sig'(1)}$ are of critical degree for the $F_i$ and
the $ F'_i$, respectively, and satisfy
$$\Res^X_F(H)=\Res^{ X'}_{ F'}( H').$$
Here $x_{\Sig(1)}$ denotes the product of the homogeneous variables $\prod_{\rho\in\Sig(1)}x_\rho$.
\end{Prop}

\begin{pf} First the sets $ U'_i=\{x\in  X' :  F'_i(x)\neq 0\}$ form 
a covering of $ X'$ since it is the pull-back of the covering $\cU$ of $X$. In particular,
the $ F'_i$ do not vanish simultaneously on $ X'$.

Now let $\Omega$ and $\Omega'$ be the Euler forms on $X$ and $ X'$, respectively. We have
$$\pi^*\big(\Omega/x_{\Sig(1)}\big)=\Omega'/x_{\Sig'(1)},$$
since both are rational extensions of the $\T$-invariant regular $n$-form 
$\frac{dt_1}{t_1}\wedge\dots\wedge\frac{dt_n}{t_n}$ on the torus, where the $t_i$ are affine coordinates.
Therefore
$$\pi^*(\omega_F(H))=\pi^*\Big(\frac{G\,\Omega/x_{\Sig(1)}}{F_0\cdots F_n}\Big)=\frac{ G'\,\Omega'/x_{\Sig'(1)}}{ F'_0\cdots F'_n}=\omega_{ F'}( H').$$
Since $\Tr_X=\Tr_{ X'}\circ\,\pi^*$ both $\omega_F(H)$ and $\omega_{ F'}( H')$ have the same toric residue. The proposition follows.
\end{pf}

\section{Residues and essential polytopes} \label{S:essential}

\begin{Def} \label{D:essential}
A collection of polytopes $P_0, \dots, P_n$ is said to be {\em
essential} if for every $I \subsetneq \{0, \dots, n\}$ the dimension
of the polytope $\sum_{i \in I} P_i$ is at least $|I|$. Given a toric
variety $X$ of dimension $n$, a collection of semi-ample degrees
$\alpha_0, \dots, \alpha_n$ is called {\it essential} if the corresponding
polytopes $P_0, \dots, P_n$ are essential.
\end{Def}

The goal of this section is to prove that the toric residue is
not identically zero if and only if the degrees $\alpha_i$
are essential.

\begin{Th} \label{T:nonzero} 
Consider degrees $\alpha_0, \dots, \alpha_n$ on a complete toric
 variety $X$. The toric residue with respect to polynomials $F_0,
 \dots, F_n$, viewed as a rational function in the coefficients of the
 $F_i$, is identically zero if and only if the $\alpha_i$ are not
 essential. For essential $\alpha_i$ there is a polynomial $H$ of
 critical degree and homogeneous of degree $1$ in the coefficients of
 each $F_i$ such that $\Res_F(H) = 1$.
\end{Th}

\begin{pf} 

The first implication is that for the non-essential degrees the toric
residue is identically 0. By \rp{birational}, we can refine $X$ to a
simplicial variety without changing the toric residue. So assume $X$
is simplicial. In this case the toric residue $\Res_F(H)$ is the sum of
the Grothendieck local residues of any $H/F_k$ with respect to the
common zeros of the remaining $F_i$ \cite[Theorem 0.4]{CCD}.

Suppose there exists a proper subset $I$ such that $\sum_{i \in I}
P_i$ has dimension less than $|I|$. Let $X_I$ be the toric variety
corresponding to $P_I=\sum_{i\in I}P_i$, and $\pi: X\to X_I$ the morphism
defined by the natural map of fans $\Sig_X\to\Sig_{P_I}$. The polynomial 
$F_i$ for $i\in I$ is the pullback of a polynomial of semi-ample degree on $X_I$
with polytope $P_i$. Clearly, generic
polynomials supported on the $P_i$, $i\in I$,
do not have a common zero on $X_I$ since $|I|>\dim X_I$.  
Thus the corresponding $\{F_i\,:i\in I\}$ do not have a common zero on $X$ for generic
coefficients. Extend $I$ to a subset of size $n$, without loss of
generality we take it to be $\{1, \dots, n\}$. For generic
coefficients $F_1, \dots, F_n$ do not have a common root. In
particular there are no local residues in the sum. So for generic coefficients
the toric residue is $0$. Since, $\Res_F$ is a rational function of the
coefficients of the $F_i$ it must be identically zero.

For the converse, the main tool is the following dual Koszul complex
of sheaves with respect to $F=(F_0, \dots, F_n)$ which appears in
numerous places including \cite{CCD}, \cite{CD}, and \cite{AK}. Given
a subset $I \subset \{0, \dots, n\}$ let $\alpha_I = \sum_{i \in I}
\alpha_i$. We have an exact sequence of sheaves
$$0 \to \cO(-\beta_0) \to \bigoplus_{i=0}^n \cO(\alpha_i - \beta_0) \to \cdots
\to \bigoplus_{|I| = p} \cO(\alpha_I - \beta_0) \to
\cdots \to \cO(\nu) \to 0,$$
where as before $\nu$ is the critical degree for $F$ and $\beta_0=\sum_{\rho}\deg(x_\rho)$.

One can take the \v{C}ech cohomology double complex and then pass to
a spectral sequence. The $E_1$ terms of this spectral sequence are:
$$E_1^{p,q} = \bigoplus_{|I| = p} H^q(X, \cO(\alpha_I - \beta_0)).$$
Because the $\alpha_i$ are essential, a result of \cite{CD} gives us:
$$E_1^{p,q} = 0 {\rm \ when\  } p+q > n, {\rm \ except\  for \ } E_1^{n+1, 0} = S_{\nu}.$$
As a consequence there is a unique top differential $d_{n+1} \ :
E_{n+1}^{0,n} \to E_{n+1}^{n+1, 0}$ which must be an isomorphism since
the spectral sequence is exact. Moreover, 
$$E_{n+1}^{0,n} = E_1^{0,n} = H^n(X, \cO(-\beta_0))$$ and
$E_{n+1}^{n+1,0}$ is a quotient of $S_{\nu}$.

So we have an induced map $S_{\nu} \to H^n(X, \cO(-\beta_0))$
which is the composition of the projection onto $E_{n+1}^{n+1,0}$
and the inverse of the isomorphism $d_{n+1}$.
We also have an isomorphism $\cO(-\beta_0)) \to \widehat{\Omega}^n$
sending a local section $s$ to $s \cdot \Omega$ where $\Omega$ is the
Euler form.  Finally there is the trace isomorphism 
$\Tr_X\,: H^n(X, \Omega^n) \to \C$. Composing all of these maps we
obtain a map $S_{\nu} \to \C$. The maps are illustrated via the diagram
below:

\begin{diagram}
S_{\nu} = E_1^{{n+1},0} &                &                       & & & &\\
\dTo                     & \rdDotsto      &                       & & & &\\  
E_{n+1}^{n+1, 0}         & \lTo^{d_{n+1}}_{\simeq} & H^n(X, \cO(-\beta_0)) & \cong & H^n(X, \widehat{\Omega}^n) & \rTo^{\Tr_X} &  \C.
\end{diagram}

We will prove that this composition is precisely the toric residue
map. In that case if we started with a differential form $\omega \in
H^n(X, \Omega^n)$ such that $\Tr_X(\omega) = 1$, it would correspond to
an element of $H^n(X, \cO(-\beta_0))$ which is mapped to an element
$h \in E_{n+1}^{n+1,0}$. Let $H$ be any element of $S_{\nu}$ lifting
$h$.  From the above constructions it would follow that $\Res_F(H) =1$ 
and the toric residue is not identically zero as desired.

Moreover, by a theorem of Weyman \cite[Chapter 3, Theorem 4.11]{GKZ},
the differential $d_{n+1}$, and therefore the element $H$ above, can
be lifted up to a (non-unique) map $H^n(X, \cO(-\beta_0)) \to S_{\nu}$
which is polynomial of degree $1$ in the coefficients of each $F_i$.

To prove that the residue map coincides with the one constructed above
we compute the cohomology terms using the \v{C}ech resolutions given
by the open cover $U_i=\{x\in X : F_i(x)\neq 0\}$. More generally,
given $J \subset \{0, \dots, n\}$ define $U_J = \cap_{j \in J} U_j$.
In this way we have the $E_0$ terms of our spectral sequence

$$E_0^{p,q} = \bigoplus_{|I| = n+1-p}\ \bigoplus_{|J|=q+1} \cO(\alpha_{\hat I}
- \beta_0)(U_J),$$ 
where $\hat I$ denotes the complement of $I$.
In terms of the cover, given a polynomial $H \in S_{\nu}$ we have
$$\frac{H}{F_0 \cdots F_n} \in \cO(-\beta_0)(U_{\{0,\dots,n\}}) = E_0^{0,n}.$$
The residue map is defined to be the trace of the
cohomology class of this latter element (after multiplying by the
Euler form). So it is enough to show that $d_{n+1}([\frac{H}{F_0\cdots
F_n}]) = [H] \in E_{n+1}^{n+1,0}$.  To compute this differential we
start with $\frac{H}{F_0 \cdots F_n} \in E_0^{0,n}$ and map it via
$d_1$ to $E_0^{1,n}$. This can be lifted via the \v{C}ech differential
$d_0$ to an element of $E_0^{1,{n-1}}$ which is further mapped to
$E_0^{2,{n-1}}$ and lifted to $E_0^{2,{n-2}}$ and so on. At the end we obtain
an element of $E_0^{n,0}$ which is mapped via $d_1$ to $E_0^{n+1,0}$.

Let $e_{IJ}$  be the basis of  $E_0^{p,q}$ and $F_I =  \prod_{i \in I}
F_i$. We have the following lemma.

\begin{Lemma} \label{L:lift} In the above mapping and lifting process, a valid
choice for the element in $E_0^{n-p,p}$ is $\sum_{|I| = p+1}\frac{H}{F_I} e_{II}$.
\end{Lemma}

\begin{pf} The base case $p = n$ is our starting element. For the
inductive step we need to show that $d_1(\sum_{|I| = p+1}\frac{H}{F_I})
e_{II} =d_0(\sum_{|I'| = p}\frac{H}{F_{I'}} e_{I'I'})$.

However, by the definitions of the Koszul and \v{C}ech morphisms it is
easy to see that both of the above elements are:
$$ \sum_{I = \{i_0, \dots, i_p\}} \sum_{j=0}^p (-1)^j\frac{H}{F_{I_j}} e_{I_jI},$$
where $I_j = I \setminus \{i_j\}$.
\end{pf}

Therefore we get the element $\sum_{i=0}^n \frac{H}{F_i} e_{ii} \in
E_0^{n,0}$.  The final Koszul differential is multiplication by $F_i$
in each factor so we are left with $(H,H,\dots, H) \in \sum_{i=0}^n
\cO(\nu)(U_i)$ which corresponds to the global section $H \in H^0(X,
\cO(\nu))$. So $H$ is a valid lifting of the image of the class of
$\frac{H}{F_0 \cdots F_n}$ under $d_{n+1}$ completing the proof.
\end{pf}

\section{Facet coloring and toric residues for monomials}\label{S:cdeg}

The theorem from the previous section guaranteed the existence of an
element of toric residue one but was completely nonconstructive.  This
section and the next one provide the framework for an explicit
combinatorial construction of such elements.  Here we recall the
definition of the facet coloring of a polytope and the relation
between the combinatorial degree of a facet coloring and the toric
residue for monomial ideals. We will also obtain the Generalized
Global Transformation Law that allows us to reduce the computation of
the toric residue for semi-ample degrees to the monomial case.

\subsection{Facet coloring}\label{S:facetcolor}
Consider an $n$-dimensional polytope $P$ in $\R^n$.
We let $\p P$ denote the boundary of $P$ and $\cF(\p P)$ the partially
ordered set (poset) by inclusion of all proper faces of $P$. We also let $2^{[n+1]}$ denote
the set of all subsets of $[n+1]=\{0,\dots,n\}$. It will be convenient for
us to equip $2^{[n+1]}$ with the inverse partial order $<$, i.e. $J<J'$ if and only if $J\supset J'$
for $J,J'\in 2^{[n+1]}$.

\begin{Def} A map of posets $C:(\cF(\p P),\subset)\to (2^{[n+1]},<)$
is called a {\it coloring of $P$ into $n+1$ colors} (or simply {\it coloring}).
The image $C(\G)$ is called the {\it set of colors} of a face $\G\in\cF(\p P)$.
We will also say that {\it $\G$ is colored by $C(\G)$}.
A coloring is called {\it simplicial} if every face $\G\in\cF(\p P)$
is colored by a non-empty proper subset of $[n+1]$. 
\end{Def}

The poset $(2^{[n+1]},<)$ can be identified with the poset of faces of 
the standard $n$-simplex:
$$\D=\{y=(y_0,\dots,y_n)\in\R^{n+1}\, :\, y_0+\dots+y_n=1,\, 0\leq y_i\leq 1\}.$$
Indeed, each non-empty proper subset $\{j_1,\dots,j_k\}\subset[n+1]$ defines the codimension $k$ face of $\D$:
$$\D_{j_1\dots j_k}=\{y\in\D\, :\, y_{j_1}=\dots=y_{j_k}=0\}.$$
Therefore, any simplicial coloring is, in fact, a map of posets $C:\cF(\p P)\to\cF(\p \D)$. 

Fix orientations of $P$ and $\D$.
Given a simplicial coloring $C$ consider a continuous piecewise linear
map $f_C:\p P\to\p\D$ such that $f_C(\G)\subset C(\G)$ for any
$\G\in\cF(\p P)$. One can show that such 
a map $f_C$ always exists and the topological degree
$\deg f_C$ does not depend on the choice of $f_C$ (see \cite{Sop}). 
We call it the {\it combinatorial degree} $\cdeg(C)$ of the simplicial
coloring $C$. The combinatorial degree is alternating in the ordering of the elements 
of $[n+1]$ as every such ordering defines an orientation of the corresponding simplex $\D$.
 
We have the following property of the combinatorial degree. 
Let $C$, $C'$ be two simplicial colorings of $P$. We say that
$C'$ {\it refines} $C$ if $C'(\G)\subset C(\G)$ for any
$\G\in\cF(\p P)$.
\begin{Prop}\label{P:refinement}\cite{Sop} 
Let $C$, $C'$ be two simplicial colorings of $P$.
If $C'$ refines $C$ then $\cdeg(C)=\cdeg(C')$.
\end{Prop}


The combinatorial degree can be computed explicitly as a signed
number of certain complete flags of faces of $P$. 
To state the precise formula we will need the following definition.
Consider a complete flag $F$ of faces of $P$
$$F:\ P^0\subset P^1\subset\dots\subset P^{n-1}\subset P^n=P,\quad\dim P^j=j.$$
For every $1\leq j\leq n$ choose a vector $e_j$ that begins at $P^0$ and
points strictly inside $P^j$. Define the {\it sign} of the flag to be $\sgn F=1$
if $(e_1,\dots,e_n)$ gives a positive oriented frame for $P$, and $\sgn F=-1$ otherwise.
It is easy to see that the sign is independent of the choice of the $e_i$.

\begin{Th}\label{T:flags} Let $C$ be a simplicial coloring of an $n$-dimensional polytope $P\subset\R^n$. 
Fix any permutation $\e$ on the elements of $[n+1]$. Then the combinatorial
degree of $C$ equals the sign of $\e$ times the number of complete flags
$$P^0\subset P^1\subset\dots\subset P^{n-1}\subset P^n=P,$$
counted with signs, such that for every $1\leq k\leq n$ the face
$P^{k-1}$ is colored by $\{\e(k),\dots,\e(n)\}$.
\end{Th}
\begin{pf} This is a particular case of \cite{Sop1}, Theorem 2.2.
\end{pf}

In particular, this theorem says that the combinatorial degree is zero unless
for every $0\leq k\leq n$ there is a facet in $P$ colored by $\{k\}$, for every
$0\leq k<l\leq n$ there is a codimension two face in $P$ colored by $\{k,l\}$, and so on.

One way to define a coloring $C:\cF(\p P)\to 2^{[n+1]}$
is to give the colors to every facet of $P$ and then extend it by taking intersections, i.e.
if $\G=\bigcap_{\nu} Q_\nu$ for some facets $Q_\nu$ then $C(\G)=\bigcup_\nu{C(Q_\nu)}$ (remember we have
the reversed order in the target). The coloring obtained 
in this way is called a {\it facet coloring}. In the present paper we will only be interested in 
facet colorings.

Next, let $P$ be a polytope in $\R^n$. Then one can consider
the poset of all flags of faces of~$P$ (chains in $\cF(\p P)$). The partial order is defined as
follows: If $F=\{\G_1\subset\dots\subset\G_k\}$ and $F'=\{\G'_1\subset\dots\subset\G'_l\}$ are two flags of faces of $P$
then $F<F'$ if and only if $\{\G'_1,\dots,\G'_l\}$ is a subset of $\{\G_1,\dots,\G_k\}$. This poset
can be realized as the poset of faces of a simple polytope $\hat P$ whose normal fan is the barycentric subdivision
of the normal fan of $P$. Indeed, one can easily see that there is a 1-1 order preserving 
correspondence between codimension $k$ faces of $\hat P$ and length $k$ flags of faces of $P$. 
In particular, facets of $\hat P$ correspond to flags of faces of length one, i.e. to faces of $P$. 

Later on we will be concerned with facet colorings not of the polytope $P$ itself,
but the polytope $\hat P$ associated with it. 
By the above, to define a facet coloring of $\hat P$
we need to assign a non-empty proper subset of $[n+1]$ to every facet of $\hat P$, hence, to every face of $P$.
Therefore, any map $C:\cF(\p P)\to 2^{[n+1]}$ defines a facet coloring $\hat C:\cF(\p\hat P)\to 2^{[n+1]}$.
(We should warn the reader, however, the map $C$ may
not be a map of posets, in general.) 
Clearly, for every flag $\G_1\subset\dots\subset\G_k$ the union
$\bigcup_i C(\G_i)$ is the set of colors of the face of $\hat P$ corresponding to this flag. We thus say
that a {\it flag $\G_1\subset\dots\subset\G_k$ is colored by $\bigcup_i C(\G_i)$}.
Furthermore, $\hat C$ is simplicial if and only if for any flag $\G_1\subset\dots\subset\G_k$ the union $\bigcup_i C(\G_i)$ is proper.




\subsection{Toric residue for monomials}
Let $X$ be a projective toric variety of dimension $n$ defined
by a lattice polytope $P$, and let $\Sig$ denote the normal fan of $P$.

Consider a collection of $n+1$ (monic) monomials $z_0,\dots,z_n$
in the homogeneous coordinate ring $S=\C[x_\rho : \rho\in\Sig(1)]$ of $X$. Assume that the
product of the variables $\prod_\rho x_\rho$ divides the product of the monomials $z_0\cdots z_n$. Then
the quotient $z_0\cdots z_n/\prod_\rho x_\rho$ has
critical degree with respect to $z_0,\dots,z_n$.

On the other hand, since the variables $x_\rho$ correspond to the 
facet normals of $P$, any
collection of monomials $z=(z_0,\dots,z_n)$ with $\prod_{\rho}x_\rho\,|\,z_0\cdots z_n$
defines a facet coloring of~$P$:
$$C_z:\cF(\p P)\to 2^{[n+1]},\quad C_z(Q_\rho)=\{i\in [n+1]\, :\, x_\rho|z_i\},$$
where $Q_\rho$ is the facet of $P$ whose inner normal generates $\rho$. 
Conversely, any facet coloring $C$ of $P$ defines a collection of squarefree monomials in $S$ whose product
is divisible by the product of the variables: $z_i=\prod_{C(Q_\rho)\ni i}x_\rho$.

If $z_0,\dots,z_n$ do not vanish simultaneously on $X$ then
the corresponding coloring $C_z$ is simplicial. Indeed, if $C_z$ is not simplicial
then there is a vertex $u$ of $P$ which is colored by $\{0,\dots,n\}$, i.e. $u\in Q_0\cap\dots\cap Q_n$
for some facets $Q_i$, such that $Q_i$ contains $i$ as one of its colors. But this implies that
the corresponding point $x_u$ on $X$ lies on the irreducible divisors $D_{\rho_0},\dots,D_{\rho_n}$,
where each $D_{\rho_i}$ is a component of the zero locus of $z_i$ on $X$, a contradiction.
 
The next theorem asserts that the combinatorial degree of $C_z$ equals
the toric residue of the quotient $z_0\cdots z_n/\prod_\rho x_\rho$.

\begin{Th}\label{T:monomres}\cite{Sop} 
Let $X$ be an $n$-dimensional projective toric
variety defined by a lattice polytope $P$. Let $z_0,\dots,z_n$ be
monomials in the homogeneous coordinate ring
$S$ such that
\begin{enumerate}
\item $\prod_\rho x_\rho\,|\,z_0\cdots z_n$,
\item $z_0,\dots,z_n$ do not vanish simultaneously on $X$.
\end{enumerate}
Then
\begin{equation}\nonumber
\Res_z({z_0\cdots z_n}/{\prod_{\rho}x_\rho})=\cdeg(C_z)
\end{equation}
where $C_z$ is the simplicial coloring of $P$ defined by $z_0,\dots,z_n$.
\end{Th}

\subsection{Reduction to toric residue for monomials}

To reduce the computation of the toric residue for arbitrary polynomials
to the case of monomials we will need the following generalized version of the
Global Transformation Law \cite{CCD}.

\begin{Th}\label{T:GTL}
Let $F_j \in S_{\alpha_j}$ and $G_j \in S_{\beta_j}$ for
$0\leq j\leq n$. Suppose
$$\sum_{j=0}^n B_{ij}F_j=\sum_{j=0}^n A_{ij}G_j,\quad 0\leq i\leq n,$$
where $B_{ij}$ and $A_{ij}$ are homogeneous of degree $\gamma_i - \alpha_j$
and $\gamma_i - \beta_j$ respectively for some fixed degrees
$\gamma_0, \dots, \gamma_n$. Assume that neither $F_0, \dots, F_n$ nor
$G_0, \dots, G_n$ vanish simultaneously on~$X$.
Let $\alpha = \sum_i \alpha_i$, $\beta = \sum_i \beta_i$,
$\gamma = \sum_i \gamma_i$, and $\nu_0=\sum_\rho\deg(x_\rho)$.
Then for any~$H\in S_{\alpha+\beta-\gamma-\nu_0}$, the polynomials
$H\det A$ and $H\det B$ are of critical degree for $F$ and $G$
respectively, and
\begin{equation}\label{e:law}
{\Res}_F(H\det A) = {\Res}_G(H\det B).
\end{equation}
\end{Th}

\begin{pf} For any $H\in S_{\alpha+\beta-\gamma-\nu_0}$
the degree of $H\det A$ is $\alpha-\nu_0$, which is the critical
degree for $F_0,\dots,F_n$. Consider the $n+1$ homogeneous polynomials
$K_i=\sum_{j=0}^n B_{ij}F_j$.
According to the Global Transformation Law (Theorem 0.1, \cite{CCD})
$$\Res_K((H\det A)\det B)=\Res_F(H\det A).$$
On the other hand, $K_i=\sum_{j=0}^n A_{ij}G_j$ and $H\det B$ has
critical degree for $G_0,\dots,G_n$. Therefore,
$$\Res_K((H\det B)\det A)=\Res_G(H\det B),$$
again by the Global Transformation Law. The theorem follows.
\end{pf}

Our reduction is then based on the following assertion.

\begin{Cor}\label{C:reduction} 
Let $X$ be an $n$-dimensional projective toric variety.
Let $F_j \in S_{\alpha_j}$ be homogeneous polynomials
not vanishing simultaneously on $X$. Suppose $y_0,\dots,y_n$ and $z_0,\dots,z_n$ are squarefree monomials
such that 
\begin{enumerate}
\item $y_0\cdots y_n={z_0\cdots z_n}/{\prod_\rho x_\rho}$,
\item $y_iF_i=\sum_{j=0}^nA_{ij}z_j$ for some $A_{ij}\in S_{\alpha_i+\deg(y_i)-\deg(z_j)}$, $0\leq i\leq n$,
\item $z_0,\dots,z_n$ do not vanish simultaneously on $X$.
\end{enumerate}
Then we have
\begin{equation}\label{e:reduction}
\Res_{F}(\det A)=\Res_{z}(y_0\cdots y_n)=\cdeg(C_z),
\end{equation}
where $C_z$ is the simplicial facet coloring defined by $z_0,\dots,z_n$.
\end{Cor}

\begin{pf} The first statement in \re{reduction} follows from \rt{GTL} and
the second statement follows from \rt{monomres}.
\end{pf}

\section{Partition matrix for polytopes and residue matrix} \label{S:partition}

\subsection{Partition matrix}
Let $P$ be a lattice polytope in $\R^n$. Consider any partition of the set of vertices of $P$ into
$n+1$ disjoint (possibly empty) subsets:
\begin{equation}\label{e:vertpart}
\Vertex(P)=V_0\sqcup\dots\sqcup V_n.
\end{equation}
Extend this partition to a partition of the set of lattice points of $P$ by adding to $V_i$
lattice points in the relative interior of faces containing a vertex from $V_i$: 
\begin{equation}\label{e:indpart}
P\cap\Z^n=M_0\sqcup\dots\sqcup M_n.
\end{equation}
Any such extension \re{indpart} will be called an {\it induced partition} of $P\cap\Z^n$ defined by the
{\it vertex partition} \re{vertpart}.

Now consider $n+1$ lattice polytopes $P_0,\dots,P_n$ in $\R^n$. For each polytope $P_i$ fix an (ordered) vertex partition
$$\Vertex(P_i)=V_{i0}\sqcup\dots\sqcup V_{in}.$$
We say that these partitions are {\it compatible} if for any permutation $\e$ of $\{0,\dots,n\}$
\begin{equation}\label{e:compatibility}
\sum_{i=0}^nV_{\e(i)i}\subset \inter\Big(\sum_{i=0}^nP_i\Big),
\end{equation}
where $\inter(P)$ denotes the relative interior of $P$.

\begin{Def} Let $P_0,\dots,P_n$ be lattice polytopes in $\R^n$. Then subsets $M_{ij}\subset P_i\cap\Z^n$,
$0\leq i,j\leq n$, form a {\it partition matrix} for $P_0,\dots,P_n$ if
$$P_i\cap\Z^n=M_{i0}\sqcup\dots\sqcup M_{in},\quad 0\leq i\leq n$$
is a collection of induced partitions defined by a compatible collection
of vertex partitions
$$\Vertex(P_i)=V_{i0}\sqcup\dots\sqcup V_{in},\quad 0\leq i\leq n.$$
\end{Def}

\begin{Rem} It is not hard to see that the compatibility condition on the $V_{ij}$ \re{compatibility}
implies the same condition on any induced partitions:
\begin{equation}\label{e:compatibility2}
\sum_{i=0}^nM_{\e(i)i}\subset \inter\Big(\sum_{i=0}^nP_i\Big).
\end{equation}
\end{Rem} 

\begin{Example}\label{Ex:partition} 
Consider three polygons $P_0$, $P_1$ and $P_2$ in \rf{partition}. We partition
their lattice points in accordance with the labels: the set
$M_{ij}$ consists points of $P_i$ labeled with $j$ ($0\leq i,j\leq 2$).

\begin{figure}[h]
\centerline{
\scalebox{0.60}{
\input{partition.pstex_t}}}           
\caption{}
\label{F:partition}
\end{figure}
Clearly these are induced partitions. To show that they are compatible it is enough
to check that for any linear functional $v\neq 0$ any three vertices 
$u_0$, $u_1$ and $u_2$ that minimize $v$ on $P_0$, $P_1$ and $P_2$, respectively, 
will not have all different labels.


\end{Example}

\subsection{Coloring matrices}
Let $P$ be a polytope in $\R^n$. Recall that every vector $v$ in the dual space $(\R^n)^*$ defines
a face $P^v$ of $P$ on which $v$ restricted to $P$ attains its minimal value.

\begin{Def} Let $M$  be a partition matrix for $P_0,\dots,P_n$. Define a map from $(\R^n)^*$ to
the set of $(0,1)$-matrices of dimension $(n+1)\times(n+1)$:
$$\cM:(\R^n)^*\to\Mat(n+1,\{0,1\})$$
where the value of $\cM$ at $v\in(\R^n)^*$ is the matrix $M^v$ whose $(i,j)$-th entry is
$$M^v_{ij}=\begin{cases}1, & \text{ if } M_{ij}\cap P_i^v\neq\emptyset,\\ 0, &\text{ otherwise}.\end{cases}$$
The matrix $M^v$ is called the {\it coloring matrix} of $v$.
\end{Def}

Informally speaking, the coloring matrix $M^v$ ``encodes'' the partitions of the lattice points of the $P_i$ 
restricted to the corresponding faces $P_i^v$. 

The compatibility condition implies that for any non-zero $v$ the coloring matrix
$M^v$ has permanent zero. Indeed, if the permanent is non-zero then there exists a permutation
$\e$ of $\{0,\dots,n\}$ such that $M^v_{\e(i)i}=1$ for all $0\leq i\leq n$. By the definition of $M^v$
this implies that for each $i$ there is a point $u_i$ in $M_{\e(i)i}$ that lies on the face $P^v_i$.
But then the sum $u_0+\dots+u_n$ gives a point on the face $P^v$ of the Minkowski sum $P=\sum_iP_i$,
which contradicts the compatibility condition \re{compatibility2}.

The following statement is known as the Frobenius-K\"onig Theorem (it is 
also equivalent to Hall's Marriage Theorem \cite{A}).

\begin{Th}\label{T:F-K}
Let $A$ be a $(0,1)$-matrix of dimension $n\times n$ with
zero permanent. Then $A$ has a submatrix of zeroes of dimension $r\times s$
for some positive $r,s$ such that $r+s=n+1$.
\end{Th}

By the above theorem for every non-zero $v$ the $(n+1)\times(n+1)$ matrix $M^v$
has a zero submatrix (not unique, in general) of dimension $r\times s$ with $r+s=n+2$.
The rows (resp. columns) of the submatrix are indexed by a subset
of $\{0,\dots,n\}$ which we denote by $I^v$ (resp. $J^v$). We thus have $|I^v|+|J^v|=n+2$
for all non-zero $v$.

Now consider a polytope $P$ whose normal fan $\Sig$ is a common refinements of the normal fans of $P_0,\dots,P_n$.
Clearly, $M^v$ is the same for all $v$ in the intersection of the cones of the $P^v_i$. Therefore, $\cM$
is constant on the cones of $\Sig$. Since cones of $\Sig$ correspond to faces of $P$ we arrive at the following definition.

\begin{Def} Let $M$ be a partition matrix for $P_0,\dots,P_n$. Let $P$ be a polytope whose normal fan is a
common refinements of the normal fans of $P_0,\dots,P_n$. Given a face $\G$ of $P$ define its {\it coloring matrix $M^\G$}
to be the coloring matrix of any $v\in\sig_\G$, where $\sig_\G$ is the cone of $\G$.
\end{Def}

We will need the following simple observation. Let $\G_1,\G_2$ be faces of $P$. Then
\begin{equation}\label{e:reverse}
\text{if }\ \ \G_1\subset \G_2\ \ \text{ then }\ \ (M^{\G_2}_{ij}=0)\Rightarrow (M^{\G_1}_{ij}=0).
\end{equation}




\subsection{Residue from a partition matrix}

Consider a projective $n$-dimensional toric variety $X$ defined by a projective fan $\Sig$. Let 
$\alpha_0,\dots,\alpha_n$ be $n+1$ semi-ample degrees on $X$ and let $P_0,\dots,P_n$ be their polytopes.

\begin{Def}\label{D:resmatrix}
Consider a collection of $n+1$ homogeneous polynomials $F=(F_0,\dots,F_n)$
of degrees $\alpha_0,\dots,\alpha_n$:
$$F_i=\sum_{u\in P_i\cap\Z^n}c_u\chi^u,\quad F_i\in S_{\alpha_i},\quad 0\leq i\leq n.$$
Given a partition matrix $M$ for $P_0,\dots,P_n$ define the {\it residue matrix} $M_F$ of $F$
to be the matrix whose entries are the homogeneous polynomials
$$F_{ij}=\sum_{u\in M_{ij}}c_u\chi^u,\quad F_{ij}\in S_{\alpha_i},\quad 0\leq i,j\leq n.$$
\end{Def}

The determinant $\det(M_F)$ is a homogeneous polynomial of degree $\alpha=\alpha_0+\dots+\alpha_n$.
Since the $\alpha_i$ are semi-ample, $\alpha$ is also semi-ample and its polytope is the Minkowski sum
$\sum_i P_i$. As follows from the definition of homogeneous coordinates (see \re{homocoord}) 
a monomial $\chi^u$ of degree $\alpha$ is divisible by all the variables if and only if the corresponding lattice
point $u$ lies in the interior of the polytope of $\alpha$. Therefore, by the compatibility condition \re{compatibility2}
every monomial in $\det(M_F)$ is divisible by all the variables, and hence the quotient
$\det(M_F)/\prod_{\rho}x_\rho$ is a homogeneous polynomial of critical degree $\alpha-\sum_\rho\deg(x_\rho)$.

\begin{Prop}\label{P:resmatrix} 
Let $\alpha_0,\dots,\alpha_n$ be semi-ample degrees on $X$ with polytopes $P_0,\dots,P_n$.
Fix a partition matrix $M$ for $P_0,\dots,P_n$. For every coloring matrix $M^{\rho}$, $\rho\in\Sig(1)$,
make any choice of an $r\times s$ zero  submatrix with $r+s=n+2$ and let its rows and columns
be indexed by subsets $I^{\rho}$ and $J^{\rho}$ of $\{0,\dots,n\}$, respectively.
Define squarefree monomials
$$y_i=\prod_{I^{\rho}\not\ni i}x_\rho,\quad z_j=\prod_{J^{\rho}\ni j}x_\rho,\quad 0\leq i,j\leq n.$$
Then for any homogeneous polynomials $F_0,\dots,F_n$ of degrees $\alpha_0,\dots,\alpha_n$
\begin{enumerate}
\item $y_0\cdots y_n={z_0\cdots z_n}/{\prod_\rho x_\rho}$,
\item $y_iF_i=\sum_{j=0}^nA_{ij}z_j$ for some $A_{ij}\in S_{\alpha_i+\deg(y_i)-\deg(z_j)}$, $0\leq i\leq n$.
\end{enumerate}
Moreover, $A_{ij}$ can be chosen so that
$$\det(M_F)/\prod_{\rho}x_\rho=\det(A),$$
where $M_F$ is the residue matrix defined by the partition matrix $M$.
\end{Prop}

\begin{pf} (1) For every $\rho\in\Sig(1)$ the variable $x_\rho$ appears in the product $z_0\cdots z_n$
with multiplicity $|J^\rho|$ and in $y_0\cdots y_n$ with multiplicity $n+1-|I^\rho|=|J^\rho|-1$
since $|I^\rho|+|J^\rho|=n+2$.

(2) For every $0\leq i\leq n$ we have 
$$F_i=\sum_{u\in P_i\cap\Z^n}c_u\chi^u.$$
We need to show that every monomial $y_i\chi^u$ is
divisible by at least one of $z_0,\dots,z_n$. Since every monomial $\chi^u$ is divisible by a vertex monomial
we can assume that $u$ is a vertex of $P_i$. Recall that in homogeneous coordinates
$$\chi^u=\prod_{\rho}x_\rho^{\langle u,v_\rho\rangle+a_\rho},$$
where $D=\sum_\rho a_\rho D_\rho$ is a representative of $\alpha_i$ (see \re{homocoord}). 
Therefore $x_\rho$ divides $\chi^u$ if and only if $\rho\not\in\sig_u$,
where $\sig_u$ is the cone of $\Sig_i$ corresponding to $u$.

The vertex $u$ is contained in $M_{ij}$ for some $0\leq j\leq n$. We show that
$z_j$ divides $y_i\chi^u$. Indeed, take any $x_\rho$ with $J^\rho$ containing $j$. 
If $i\in I^\rho$ then $M^\rho_{ij}=0$. From the definition of $M^\rho$ it follows that $P_i^{v_\rho}$ does not
contain the vertex $u$, i.e. $\rho\not\in\sig_u$ and so $x_\rho|\chi^u$ by above. If
$i\not\in I^\rho$ then $x_\rho|y_i$ by the definition of $y_i$. 

The above argument shows that $y_iF_{ij}=A_{ij}z_j$ for some homogeneous polynomial $A_{ij}$. Taking the determinant
we obtain $y_0\cdots y_n\det(M_F)=z_0\cdots z_n\det(A)$. Now the last statement follows from part (1).
\end{pf}

The above proposition
shows that given a partition matrix $M$, any choice of zero submatrices in $M^\rho$,
for $\rho\in\Sig(1)$, defines a collection of squarefree monomials $y_0,\dots,y_n$ and $z_0,\dots,z_n$
that satisfy the conditions (1) and (2) of \rc{reduction}.
If the facet coloring $C_z$ defined by the monomials $z_0,\dots,z_n$ is simplicial the condition (3) of \rc{reduction} is satisfied
and that would imply the result of \rt{main1}, namely that the residue of $M_F$ equals the combinatorial 
degree of $C_z$. However, there are examples of $P_0,\dots,P_n$ when the condition (3) fails
no matter how one chooses a partition matrix and zero submatrices. To avoid this obstruction
we are going to change the variety $X$ by taking the barycentric refinement of its fan $\hat\Sig\to\Sig$.
This gives a birational morphism $\hat X\to X$ which allows us to transfer our construction to the variety $\hat X$
(see \rp{birational}). The advantage of this is that for any partition matrix $M$ there is a canonical choice of a zero submatrix in every
coloring matrix $M^\rho$, for $\rho\in\hat\Sig$, which guarantees that the corresponding monomials
$\hat z_0,\dots,\hat z_n$ do not vanish simultaneously on $\hat X$.

\section{Canonical colorings}\label{S:canoncolor}


Let $P_0,\dots,P_n$ be $n+1$ lattice polytopes in $\R^n$ and $\Sig_0,\dots,\Sig_n$ their
normal fans. Let $P$ be any polytope whose normal fan $\Sig$ is a common refinement of the $\Sig_i$.
Given a partition matrix we will define a canonical facet coloring of
a polytope $\hat P$ whose normal fan $\hat\Sig$ is the barycentric refinement of $\Sig$.
We will then prove that this coloring is simplicial. 
This will allow us to define monomials $\hat z_0,\dots,\hat z_n$
on the toric variety corresponding to $\hat P$ that satisfy all the conditions of \rc{reduction} and
thus obtain our main result (\rt{main1}).

\subsection{Canonical coloring}\label{S:canonical} 

Let $M$ be a partition matrix for polytopes $P_0,\dots,P_n$ and
let the polytopes $P$ and $\hat P$ be as above. As mentioned in \rs{facetcolor}, to define a facet coloring of $\hat P$ it
suffices to assign a subset $C(\G)\subset \{0,\dots,n\}$ to every face $\G$ of $P$. We will start by describing
all possible candidates for $C(\G)$, so called admissible colorings of $\G$.

\begin{Def} Let $\G$ be a face of $P$ and $M^\G$ its coloring matrix. A subset 
$J\subset\{0,\dots,n\}$ is called an {\it admissible coloring of $\G$} if $M^\G$
contains an $r\times s$ zero submatrix with $r+s=n+2$ whose columns are indexed by $J$.
\end{Def}

It turns out that the set of admissible colorings of a face possesses very nice properties.

First, for any flag of faces of $P$ we have the reversed inclusion of the corresponding sets of admissible colorings:
\begin{equation}\label{e:reverse2} 
\text{ If }\ \  \G_1\subset\dots\subset \G_k\ \ \text{ then }\ \ \cJ_1\supset\dots\supset\cJ_k,
\end{equation}
where $\cJ_i$ is the set of admissible colorings of $\G_i$. 
Indeed, \re{reverse} implies that every zero submatrix in the coloring matrix of $\G_i$ is also a zero submarix
of the coloring matrix of $\G_{i-1}$.

Second, let $M^\G$ be the coloring matrix of a face $\G\subset P$. (In what follows we will only
use that $M^\G$ is an $(n+1)\times(n+1)$ matrix with $(0,1)$-entries and zero permanent.)
Denote by $\cB$ the set of all zero submatrices $B$ in $M^\G$
of dimension $r\times s$ such that $r+s=n+2$. This set is non-empty by \rt{F-K}.
For $B\in\cB$ we let $I(B)$ (resp. $J(B)$) denote the subset in $\{0,\dots,n\}$
of indices of rows (resp. columns) of $B$. We have the following lemma.

\begin{Lemma}\label{L:lattice}
Let $B_1, B_2\in\cB$. Then there is $B\in\cB$ such that either $J(B)=J(B_1)\cup J(B_2)$
or $J(B)=J(B_1)\cap J(B_2)$.
\end{Lemma}
\begin{pf} Let $B'$ be the submatrix whose rows are indexed by $I(B_1)\cap I(B_2)$ and
whose columns are indexed by $J(B_1)\cup J(B_2)$. Clearly $B'$ is a zero submatrix. Similarly,
let $B''$ be the zero submatrix with rows indexed by $I(B_1)\cup I(B_2)$ and
columns indexed by $J(B_1)\cap J(B_2)$. Denote $r_i=|I(B_i)|$, $r_\cap=|I(B_1)\cap I(B_2)|$,
and $r_\cup=|I(B_1)\cup I(B_2)|$. By the inclusion/exclusion formula $r_\cup+r_\cap=r_1+r_2$. Similarly,
$s_\cup+s_\cap=s_1+s_2$, where $s_i=|J(B_i)|$, $s_\cap=|J(B_1)\cap J(B_2)|$, and $s_\cup=|J(B_1)\cup J(B_2)|$.
Summing up these two equations we obtain $$(r_\cap+s_\cup)+(r_\cup+s_\cap)=(r_1+s_1)+(r_2+s_2)=2(n+2).$$
Therefore, either $r_\cap+s_\cup\geq n+2$ or $r_\cup+s_\cap\geq n+2$. In other words, either $B'$ or $B''$
contains a zero submatrix $B$ with $r+s=n+2$, as required.
\end{pf}

\begin{Rem} The above lemma means that if $J_1$ and $J_2$ are two admissible colorings of $\G$ then
either $J_1\cap J_2$ or $J_1\cup J_2$ is also an admissible coloring. As follows from the proof,
a slightly stronger statement is true: If $J_1\cup J_2$ is not an admissible coloring then
any single color can be removed from $J_1\cap J_2$ and the remaining set
will still be an admissible coloring of~$\G$.
\end{Rem}

\begin{Lemma} \label{L:canonical}
Let $\cB$ be as above and consider the partially ordered  by inclusion set
$$\cJ=\{J(B)\subset\{0,\dots,n\} : B\in\cB\}.$$
Let $\cJ_\cup$ be the set of maximal elements, and $\cJ_\cap$ the set of minimal elements of $\cJ$. 
Then the subsets 
$$c=\bigcup_{J\in\cJ_\cap}J\quad\text{and}\quad C=\bigcap_{J\in\cJ_\cup}J$$
belong to $\cJ$ and satisfy $c\subset C$.
\end{Lemma}

\begin{pf} 
To prove $C\in\cJ$ we 
show that $J_1\cap\dots\cap J_k\in\cJ$ for any $J_i\in\cJ_\cup$, $1\leq i\leq k$. We proceed by induction.
The case $k=1$ is trivial. Assume $J=J_1\cap\dots\cap J_k\in\cJ$ and let $J_{k+1}\in\cJ_\cup$.
If $J\cap J_{k+1}\in\cJ$ we are done, otherwise $J\cup J_{k+1}\in\cJ$ by \rl{lattice}.
Since $J_{k+1}$ is maximal we have $J\subset J_{k+1}$, i.e. $J=J\cap J_{k+1}=J_1\cap\dots\cap J_k\cap J_{k+1}\in\cJ$.
Similar arguments show that $c\in\cJ$. 

To show $c\subset C$ it is enough to notice that for any $J\in\cJ_\cap$ and any $J'\in\cJ_\cup$ we have $J\subset J'$. 
Indeed, either $J\cap J'\in\cJ$ and so $J=J\cap J'\subset J'$ by minimality of $J$, 
or $J\cup J'\in\cJ$ and so $J\subset J\cup J'=J'$ by maximality of $J'$.
\end{pf}

The above lemma supplies us with two canonical coloring of a face $\G$:

\begin{Def} Let $M^\G$ be the coloring matrix of a face $\G\subset P$ and $\cJ_\G$ the set of all 
admissible colorings of $\G$.
Maximal (minimal) elements of $\cJ_\G$ are called {\it maximal (minimal) colorings} of $\G$.
The union $c(\G)$ of minimal colorings is called the {\it minimal canonical coloring} of $\G$. The intersection
$C(\G)$ of maximal colorings is called the {\it maximal canonical coloring} of $\G$. 
\end{Def}

\begin{Example} Consider the three polygons $P_0$, $P_1$ and $P_2$ from \rex{partition}.
Let $\G$ be the horizontal edge of the Minkowski sum $P=P_0+P_1+P_2$. Then it has the coloring 
matrix
$$M^\G=\left[\begin{matrix}0 & 1 & 0 \\ 1 & 1 & 0 \\ 1 & 0 & 0\end{matrix}\right].$$
We get $\cJ_\G=\left\{\{2\}\right\}$ and $c(\G)=C(\G)=\{2\}$. Next let $\G'$ be the edge of $P$ with $45^\circ$
slope. (It is the sum of the highest vertex of $P_0$ and the two edges of $P_1$ and $P_2$ of slope $45^\circ$.)
Its coloring matrix is
$$M^{\G'}=\left[\begin{matrix}0 & 0 & 1 \\ 1 & 0 & 1 \\ 0 & 0 & 1\end{matrix}\right].$$
This time we have $\cJ_{\G'}=\left\{\{1\},\{0,1\}\right\}$ and $c(\G')=\{1\}$, $C(\G')=\{0,1\}$.

To obtain less trivial example we need to consider the case $n=3$. Here is an example of a coloring
matrix whose set of admissible colorings has more than one maximal (and minimal) element.
$$M^\G=\left[\begin{matrix}1& 0 & 0 & 0 \\ 0 & 0 & 1 & 0 \\ 0 & 0 & 1 & 0 \\ 1 & 0 & 0 & 0\end{matrix}\right].$$
Indeed, the set of admissible colorings is $\cJ=\left\{\{1\},\{3\},\{1,3\},\{0,1,3\},\{1,2,3\}\right\}$.
Therefore, $c=\{1\}\cup\{3\}=\{1,3\}$ and $C=\{0,1,3\}\cap\{1,2,3\}=\{1,3\}$.
\end{Example}

According to the discussion in \rs{facetcolor} the maps $c:\G\mapsto c(\G)$ and $C:\G\mapsto C(\G)$ defined above give rise to two 
facet coloring $\hat c$ and $\hat C$ of $\hat P$ which we call the minimal and maximal canonical facet
colorings of $\hat P$, respectively. 

It is easy to see that under $\hat C$ no facet of $\hat P$ gets all the colors. Indeed, for any face $\G\subset P$
its coloring matrix  $M^\G$ cannot contain zero rows, thus every maximal coloring $J\in\cJ_{\G}$ is a 
proper subset of $\{0,\dots,n\}$.
The next theorem shows that an even stronger statement is true: no face of $\hat P$ gets all the colors.

\begin{Th}\label{T:simplicial} 
The maximal and minimal canonical facet colorings of $\hat P$ are simplicial and have the
same combinatorial degree.
\end{Th}

\begin{pf} Recall from \rs{facetcolor} that to prove $\hat C$ simplicial we need to show that for any maximal flag of faces of $P$
$$\G_0\subset\dots\subset \G_{n-1},\quad \dim\G_i=i$$
the union $\bigcup_{i=0}^{n-1}C(\G_i)$ is a proper subset of $\{0,\dots,n\}$. We will
prove by induction that for any $n-1\geq k\geq 0$
\begin{equation}\label{e:ih}
\bigcup_{i=k}^{n-1}C(\G_i)\subset J_k,
\end{equation}
for some maximal coloring $J_k$ of $\G_k$. For $k=0$ this 
implies the statement of the theorem, since $J_0$ is a proper subset of $\{0,\dots,n\}$.
The base $k=n-1$ is clear since $C(\G_{n-1})$ is the intersection of maximal elements of $\cJ_{n-1}$.
For the inductive step assume that \re{ih} is true
for some maximal $J_k\in\cJ_k$. By \re{reverse2} $\cJ_{k}\subset\cJ_{k-1}$, thus there exists a 
maximal element $J_{k-1}\in\cJ_{k-1}$
such that $J_k\subset J_{k-1}$. Also $C(\G_{k-1})\subset J_{k-1}$, by definition. This together with \re{ih} gives
$$\bigcup_{i=k-1}^{n-1}C(\G_i)\subset J_{k-1},$$
as required. 

By \rl{canonical} $c(\G)\subset C(\G)$ for any face $\G\subset P$. Therefore, $\hat c$ is also simplicial.
Finally, $\cdeg(c)=\cdeg(C)$ follows from \rp{refinement}.
\end{pf}

When the maximal canonical coloring of a face consists of a single element we can say more
about admissible colorings of this face:

\begin{Lemma}\label{L:singlecolor} Suppose a face $\G\subset P$ is 
maximally canonically colored by a single color $C(\G)=\{k\}$.  
Then this is the only admissible coloring of $\G$. Moreover,
any face containing $\G$ is also singly canonically colored by $\{k\}$ while every
subface of $\G$ is canonically colored by a set containing $k$.
\end{Lemma}

\begin{pf} Suppose $J_1,\dots,J_s$, $s\geq 2$, are maximal colorings of $\G$ such that
$\{k\}=J_1\cap\dots\cap J_s$, but $\{k\}\subsetneq J_1\cap\dots\cap J_{s-1}$.
By the proof of \rl{canonical} $J=J_1\cap\dots\cap J_{s-1}$ is an admissible coloring
of $\G$. By the remark after \rl{lattice} $\G$ can either be colored by $J\cup J_s$
or else by $J\cap J_s$ with any single color removed. The first is a coloring strictly
larger than $J_s$ which is impossible since $J_s$ is maximal. The second is empty since
$J\cap J_s$ is already a single color. Both are contradictions. Thus, the unique
maximal coloring of $\G$ is $\{k\}$ which is therefore the only admissible coloring.

If $\G\subset\G'$ then $\cJ_\G\supset\cJ_{\G'}$ by \re{reverse2} and, hence, $\cJ_{\G'}=\cJ_\G=\{k\}$.
If $\G\supset\G''$ then $\{k\}$ is an admissible coloring of $\G''$. But $\{k\}$ can
be appended to any coloring of $\G''$. Thus $k$ is contained in every maximal coloring of $\G''$,
i.e. in the maximal canonical coloring of $\G''$.
\end{pf}

\subsection{Main theorem}

We now turn back to residues and prove the result of \rt{main1}. As before $X$ is a complete $n$-dimensional
toric variety defined by a fan $\Sig$ and $\alpha_0,\dots,\alpha_n$ are semi-ample
degrees with polytopes $P_0,\dots,P_n$. We can assume that $X$ is projective
and take $P$ to be the polytope of an ample divisor on $X$ (If $X$ is not projective
it can be dominated birationally by a projective toric variety. This will not affect
the toric residue computation by \rp{birational}.) We also let $\hat P$ denote a polytope
whose normal fan is the barycentric subdivision of $\Sig$.

Let $M$ be a partition matrix for $P_0,\dots,P_n$. According to \rs{canonical} $M$ produces
a map $C:\cF(\p P)\to 2^{[n+1]}$ which assigns to every proper face $\G$ of $P$ its 
maximal canonical coloring $C(\G)$. The induced canonical facet coloring $\hat C$ of $\hat P$
is simplicial by \rt{simplicial}. The next theorem says that for any $F_0,\dots,F_n$
of degrees $\alpha_0,\dots,\alpha_n$ the determinant of the residue matrix $M_F$ (see \rd{resmatrix})
gives an element whose residue is the combinatorial degree of $\hat C$.

\begin{Th} Let $X$ be a complete toric variety of dimension $n$. Let $\alpha_0,\dots,\alpha_n$
be semi-ample degrees and $P_0,\dots,P_n$ their polytopes. Consider a partition matrix $M$
for $P_0,\dots,P_n$. For any collection of homogeneous polynomials $F_0,\dots,F_n$ of
degrees $\alpha_0,\dots,\alpha_n$ consider the corresponding residue matrix $M_F$.
Then the residue of $\det(M_F)/\prod_{\rho}x_\rho$ is equal to the combinatorial degree 
of the canonical facet coloring of $\hat P$:
$$\Res_F\big(\det(M_F)/\prod_{\rho}x_\rho\big)=\cdeg(\hat C).$$
\end{Th}

\begin{pf} First notice that we can work on the variety $\hat X$ defined by the polytope $\hat P$.
Indeed, let $\pi:\hat X\to X$ be the birational morphism defined by the barycentric refinement $\hat\Sig\to\Sig$
and $\pi^*:S\to\hat S$ the induced homomorphism of homogeneous coordinate rings.
Then each polynomial $\hat F_i=\pi^*(F_i)$ is of semi-ample degree $\hat\alpha_i=\pi^*(\alpha_i)$ and by \rp{birational} 
$$\Res^X_F(H)=\Res^{\hat X}_{\hat F}(\hat H),$$
where $H=\det(M_F)/\prod_{\rho}x_\rho$ and $\hat H=\pi^*(\det(M_F))/\prod_{\hat\rho}x_{\hat\rho}$, for $\hat\rho\in\hat\Sig(1)$. 

Since the degrees $\hat\alpha_i$ have the same polytopes $P_i$ we did not change the partition matrix and
the pullback $\pi^*(M_F)$ is the residue matrix $M_{\hat F}$ for the $\hat F_i$. Therefore we can apply
\rp{resmatrix} for the canonical facet coloring of $\hat P$
to obtain squarefree monomials $\hat y_0,\dots,\hat y_n$ and $\hat z_0,\dots,\hat z_n$ in $\hat S$
which satisfy
\begin{enumerate}
\item $\hat y_0\cdots \hat y_n={\hat z_0\cdots \hat z_n}/{\prod_{\hat \rho} x_{\hat \rho}}$,
\item $\hat y_i\hat F_i=\sum_{j=0}^n\hat A_{ij}\hat z_j$ for some $\hat A_{ij}\in \hat S_{\alpha_i+\deg(\hat y_i)-\deg(\hat z_j)}$, $0\leq i\leq n$,
\item $\hat z_0,\dots,\hat z_n$ do not vanish simultaneously on $\hat X$,
\item $\det(M_{\hat F})/\prod_{\hat\rho}x_{\hat\rho}=\det\hat A$.
\end{enumerate}
(Part (3) follows since the $\hat z_i$ define the canonical facet coloring $\hat C$ of $\hat P$ which is simplicial according to \rt{simplicial}.)
By \rc{reduction} 
$$\Res_{\hat F}(\det\hat A)=\cdeg(\hat C),$$
which completes the proof.
\end{pf}

\section{Locally Unmixed Degrees}\label{S:locallyunmixed}

In this section we consider the special case when the $n+1$ polytopes share a complete flag of faces.
An essential family of degrees with such collection of polytopes is called locally unmixed.
We show that for any family of locally unmixed degrees one can write an explicit partition matrix yielding
an element of residue $\pm 1$ (\rt{loc_unmixed}).

\begin{Def} Polytopes $P_0, \dots, P_m \subset \mathbb{R}^n$ are said
to {\em share a complete flag} if for each $P_i$ there is a complete flag of
faces:
$$P_i^0 \subset P_i^1 \subset \cdots \subset P_i^{n-1}, \quad \dim
P_i^j = j,$$such that the sums of the corresponding entries $P^j
=\sum_{i=0}^m P_i^j$ form a complete flag of faces of $P = \sum_{i=0}^m P_i$:
$$P^0 \subset P^1 \subset \cdots \subset P^{n-1} \subset P^n = P, \quad \dim P^j = j.$$
\end{Def}

An immediate consequence of the above definition is that if $I$ is any
non-empty subset of $\{0, \dots, m\}$ we can similarly define $P_I =
\sum_{i \in I} P_i$ such that the $P_I^j = \sum_{i \in I} P_i^j$
also form a complete flag of faces of $P_I$:
$$P_I^0 \subset P_I^1 \cdots \subset P_I^{n-1}.$$

\begin{Def} Let $X$ be a complete toric variety of dimension $n$.
An essential family of semi-ample degrees $\alpha_0, \dots, \alpha_n$ is said to be {\em
locally unmixed} if the corresponding polytopes $P_0, \dots, P_n$
share a complete flag.
\end{Def}

Note that the $P_i$ themselves may be only $n-1$ dimensional, although
at least two of them must be $n$-dimensional since the family is essential.

\begin{Th}\label{T:loc_unmixed} Let $\alpha_0, \dots, \alpha_n$ be locally unmixed degrees
on $X$. Define partitions:
$$M_{ij} = \{ u \in P_i^j \cap \Z^n \quad { \rm with } \quad u \notin
P_i^{j-1} \cap \Z^n\}.$$ This is a compatible collection of partitions
and the corresponding residue matrix gives an element of residue $\pm
1$ for any homogeneous polynomials $F_i \in S_{\alpha_i}$ not vanishing simultaneously on $X$.
\end{Th}

Before we begin the proof let us illustrate the partition using the following 3-dimensional example.

\begin{Example} Consider four 3-dimensional polytopes $P_0$, $P_1$, $P_2$ and $P_3$ as in \rf{polytopes}.
They share a complete flag of faces. Indeed, each of them has a face with inner normal $(0,0,-1)$ and
this face has an edge along the vector $(1,0,0)$. Then for each $0\leq i\leq 3$ set $M_{i0}$ consists of the
vertex of the flag (point marked as ``0''), set $M_{i1}$ consists of the other lattice 
points on the edge of the flag (lattice points marked as ``1''), set $M_{i2}$ consists of the
lattice points on the face of the flag, but not on the edge (points marked as ``2''), and
the rest of the lattice points constitute $M_{i3}$ (points marked as ``3''). 

\begin{figure}[h]
 \[
   \begin{picture}(345,100)
     \put(1,1){\epsfxsize=12.3cm\epsffile{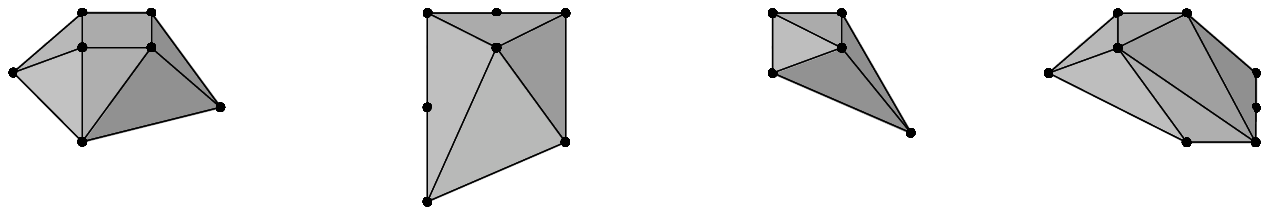}}
     \put(15,90){$P_0$} \put(115, 90){$P_1$} \put(210,90){$P_2$} \put(310, 90){$P_3$}
     \put(14,60){$0$}  \put(112, 60){$0$}   \put(209, 60){$0$} \put(305, 60){$0$}
     \put(36,60){$1$}  \put(134, 60){$1$} \put(156, 60){$1$}  \put(231, 60){$1$} \put(327, 60){$1$}
     \put(23,47){$2$}  \put(48, 47){$2$} \put(124, 41){$2$}  \put(243, 45){$2$} \put(319, 42){$2$}
     \put(-8,35){$3$}  \put(62,25){$3$} \put(108, 25){$3$} \put(205, 35){$3$}  
     \put(284, 35){$3$} \put(354, 39){$3$} \put(14,8){$3$}  \put(108, -1){$3$} 
     \put(158, 12){$3$}  \put(255, 12){$3$} \put(320, 12){$3$} \put(354, 26){$3$} \put(354, 12){$3$}   
   \end{picture}
\]
\caption{}
\label{F:polytopes}
\end{figure}

\end{Example}

We start with a simple lemma.

\begin{Lemma}\label{L:sum}
Let $P$ be a polytope of dimension $n$ and $P'$ a polytope 
such that $P+P'$ is also a polytope of dimension $n$. 
For any facet $Q$ of $P$ there is a unique face 
$\G'$ of $P'$ such that $Q+\G'$ is a proper face (in fact a facet) of $P+P'$. Hence, if $u \in P$ is a point
in the relative interior of $Q$ and $u' \in P'$ is not on
the corresponding face $\G'$ of $P'$ then $u+u'$ is in the interior of $P + P'$.
\end{Lemma}

\begin{pf}

Let $\R^n$ be the affine span of $P$. Since $P+P'$ is also
$n$-dimensional, we must have $P' \subset \mathbb{R}^n$ and $P+P' \subset
\mathbb{R}^n$. For any facet $Q$ of $P$ there is a unique linear functional
(up to scaling) $v_Q \in (\mathbb{R}^n)^*$ minimized on $Q$ in $P$.  Let
$\G'$ be the unique maximal face of $P'$ on which $v_Q$ is minimized.
The Minkowski sum $Q+\G'$ is the facet of $P+P'$ on which $v_Q$ is
minimized and conversely any face of $P+P'$ with $Q$ a summand must
minimize $v_Q$ and so must be $Q+\G'$.

For the second statement, note that $Q$ is the only face of $P$
containing $u$. By the first part, every face $\G''$ of $P'$ such that
$Q+\G''$ is contained in a proper face of $P+P'$ must have $\G'' \subset \G'$.  As a
consequence if $u'$ is not on $\G'$, hence not on any such $\G''$, $u+u'$ is
not contained in any proper face of $P+P'$.
\end{pf}

\begin{proof}[Proof of \rt{loc_unmixed}]
The lattice point partitions $M_{ij}$ are induced from the vertex
partitions obtained from the same rule restricted to the vertices of
the $P_i$. To show that $M$ is a partition matrix we must show
that:
$$\sum_{i=0}^nM_{\e(i)i} \subset \inter \Big(\sum_{i=0}^nP_i\Big)$$
for any permutation $\e$ of $\{0, \dots, n\}$.
We will show by induction that
$$\sum_{i=0}^jM_{\e(i)i} \subset \inter \Big(\sum_{i=0}^jP_{\e(i)}^j
\Big)$$ for $j = 0, \dots, n$.  The case $j=n$ is our desired
result. Let $I(j) = \{\e(0), \dots, \e(j)\} \subset \{0, \dots,n\}$.
Hence, the right hand side is $P_{I(j)}^j$, a polytope of dimension $j$.

The case $j=0$ is trivial. For the induction, we assume
$$\sum_{i=0}^{j-1}M_{\e(i)i} \subset \inter \big(P_{I(j-1)}^{j-1}\big).$$

Next, $P_{I(j-1)}^{j-1}$ is a facet of $P_{I(j-1)}^j$ (the case
$j=n$ requires that $P_{I(n-1)}^n$ is actually $n$-dimensional), so
we apply \rl{sum}. Any point in $\sum_{i=0}^{j-1}M_{\e(i)i}$ lies in
the interior of $P_{I(j-1)}^{j-1}$, and any point in $M_{\e(j)j}$
does not lie on the associated face $P_{\e(j)}^{j-1}$ of
$P_{\e(j)}^j$. Therefore, by \rl{sum}, any point in $\sum_{i=0}^j
M_{\e(i)i}$ is in the (relative) interior of $P_{I(j-1)}^j +
P_{\e(j)}^j = P_{I(j)}^j$ as desired.

To show that the combinatorial degree of the maximal canonical
coloring of $\hat{P}$ is $\pm 1$ we apply \rt{flags}.
Recall that a face of codimension $k$ of $\hat{P}$ is a flag of $k$
faces $\G_{i_1} \subset \G_{i_2} \subset \cdots \subset \G_{i_k}$ of
$P$.  We show that
there is only one complete flag of faces of $\hat{P}$ colored $(\{n\},
\{n, n-1\}, \dots, \{n, \dots, 1\})$, namely $(P^{n-1}, (P^{n-1}, P^{n-2}),
\dots, (P^{n-1}, \dots, P^0))$.

To do this we prove a few simple lemmas:

\begin{Lemma} \label{L:lem1}
The maximal canonical coloring of the face $P^j\subset P$ for $j < n$ is $\{j+1,\dots, n\}$.
\end{Lemma}

\begin{pf}
The polytope $P^j$ is the Minkowski sum of $P^j_0, \dots, P^j_n$, and
each $P^j_i$ contains precisely all of the lattice points in $M_{ik}$
for $k = 0, \dots, j$. Thus, the corresponding coloring matrix for
$P^j$ has all $1$'s in columns $0, \dots, j$ and all $0$'s in columns
$j+1, \dots, n$. It follows immediately that the only maximal coloring
is $\{j+1,\dots, n\}$, as desired.
\end{pf}

\begin{Lemma} \label{L:lem2}
The maximal canonical coloring of any proper subface of $P^j$ other than 
$P^{j-1}$ contains some color $k$ with $k < j$.
\end{Lemma}

\begin{pf}
Let $\G$ be a proper subface of $P^j$.  We decompose $\G$ as the
 Minkowski sum $\G_0 + \cdots + \G_n$ where each $\G_i$ is a subface of
 $P_i^j$.  If $j=n$ and if $\G$ were a counterexample to the lemma it
 would have to be colored just $\{n\}$.  The last column of its
 coloring matrix is $0$.  Consequently $\G_i$ contains no points of
 $M_{in}$ and so is entirely contained in $P_i^{n-1}$. So we can
 reduce to the case $j < n$ and assume that each $\G_i$ is a {\em
 proper} subface of $P_i^j$.

Now assume $\G \neq P^{j-1}$. We show that we can take $k$ to be the
 smallest number such that for all $i$, $P_i^k \nsubseteq \G_i$. If
 $P_i^{j-1} \subset \G_i$ for some $i$ then as $\G_i$ is a proper face
 of $P_i^j$ we must have $\G_i = P_i^{j-1}$.  This is a facet of
 $P_i^j$, so repeated applications of \rl{sum} show that every other
 summand $\G_{i'} = P_{i'}^{j-1}$ and so $\G = P^{j-1}$, a
 contradiction. Therefore, $k \leq j-1$.

By hypothesis, for some $i$, $P_i^{k-1} \subset \G_i$ but $P_i^k
\nsubseteq \G_i$. In particular $\G_i \cap M_{ik} = \emptyset$. If
$\G_{i'} \cap M_{i'k} \neq \emptyset$ for some $i' \neq i$, then
another application of \rl{sum} shows that $\G_i + \G_{i'}$ contains a
point in the relative interior of $P_i^k + P_{i'}^k$ and so must
contain the entire face.  But this would imply $\G_i$ contains
$P_i^k$, a contradiction.  Therefore, in the coloring matrix of $\G$
coming from $M$, the entire column $C_k$ is 0 and so $k$ is part of
the canonical maximal coloring of $\G$ as desired.
\end{pf}

Our desired  result now follows by induction.  By  \rl{lem1},
 $P^{n-1}$ is  colored just  $\{n\}$ and by \rl{lem2} it  is the
 only such face of $P$ (facet of $\hat{P}$).  Inductively, the face of
 $\hat{P}$ given by the flag of faces $(P^{n-1}, P^{n-2}, \dots, P^j)$
 in $P$ is colored $\{n, \dots,  j+1\}$. For the next step we must
 add a subface of $P^j$ to the flag with $j$ the only new color. But
 by \rl{lem1} and \rl{lem2}, the  only such  subface is $P^{j-1}$.
\end{proof}

\section{Dimension Two}\label{S:dim2}

In this section we prove that matrices whose determinant have residue $\pm 1$
can be found for almost all essential, two-dimensional families of degrees.

Recall the definition of essential in \rd{essential}. In the special case $n=2$, essential means that
no $P_i$ is zero dimensional, and while some or all of the $P_i$ may
be one dimensional line segments, no two such are parallel line
segments. We will show we can always find a residue matrix that gives an element
of residue $\pm 1$ in all but one exceptional case.

\begin{Def} Degrees $\alpha_0, \alpha_1, \alpha_2$ are {\em exceptional}
if for two of them, $\alpha_i$ and $\alpha_j$, the corresponding
polygons $P_i$ and $P_j$ are 1-dimensional, and the third $\alpha_k$
is an ample divisor on the toric variety defined by $P_i + P_j$.
\end{Def}

\begin{Th} Let $\alpha_0, \alpha_1, \alpha_2$ be an essential,
non exceptional family of degrees on a toric surface $X$. There exists
a partition matrix for the $\alpha_i$ which yields an element of
residue $\pm 1$ for every set of $F_i \in S_{\alpha_i}$ without a
common root.
\end{Th}

Note that the codimension 1 theorem for the critical degrees has been
proved by Cox and Dickenstein \cite{CD} when all $\alpha_i$ are full
dimensional. Such a case, of course, will never be exceptional.  It
is, however possible for the critical degree to be of codimension 1,
in which case the residue map is an isomorphism, and still be exceptional.
See \rex{exbad} below.

\begin{pf}
Let $P_0, P_1, P_2$ be the corresponding polygons and $P = P_0 + P_1
+ P_2$ their Minkowski sum. Every edge $e$ of $P$ is the sum of edges
from one or more of the $P_i$ and vertices from the others. Label an
edge by a subset of $\{0, 1, 2\}$ consisting of those polygons
for which the summand of $e$ is an edge. Now consider consecutive
edges of $P$. Proceed until we have a sequence containing all three
labels $0, 1, 2$. Take the smallest subsequence with this property.
We then have the following cases:

\begin{enumerate}

\item The sequence has length $1$, so there is a single edge
labeled $\boxed{0\,1\,2}$. This will be the locally unmixed case.

\item The sequence has length $2$. Up to relabeling and change of
direction the sequence will be either:
\begin{enumerate}
\item $\boxed{0\,1}\,,\, \boxed{2}$ or
\item $\boxed{0\,1}\,,\, \boxed{1\,2}$.
\end{enumerate}
Such sequences will be called {\em partially unmixed}.

\item The sequence has length $3$ or more. All such sequences can
be represented as

$\boxed{0\, \textcolor[gray]{0.4}{1}}\,,\, \boxed{1}\,,\, 
\textcolor[gray]{0.4}{\boxed{1}}\,, \dots,\, \textcolor[gray]{0.4}{\boxed{1}}\,,\, \boxed{\textcolor[gray]{0.4}{1}\, 2}$.

The numbers in gray may or may not occur. That is to say, the
first term must contain label $0$ but may or may not contain
$1$. Similarly for the final term must contain label $2$ and possibly
also label $1$.  There is at least one term labeled just $1$ in the
middle, but there may be others. Altogether, sequences of this type
will be called {\em generically mixed}.
\end{enumerate}

\begin{description}
\item[Case 1] The three degrees share an edge, hence share a complete flag
consisting of this edge and either of the two vertices. So the
polygons are locally unmixed in the sense of the previous
section. Therefore, we know we can always find a partition yielding a
residue $1$ matrix. We illustrate the partition via the
diagram in \rf{lup}.

\begin{figure}[h]
\centerline{
\scalebox{0.60}{
\input{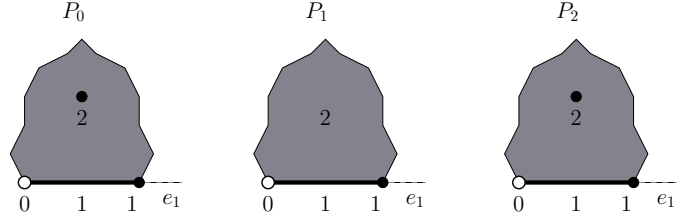}
             }
           }
  \caption{Case 1: Locally Unmixed Partition}
         \label{F:lup}
     \end{figure}

Each of the three figures represents one of the three polygons.  The
edge $e_1$ on the bottom is shared by all three polygons.  The
white vertex in each $P_i$ is in partition set $M_{i0}$ hence marked
as ``$0$'' in the diagram, the rest of the lattice points on $e_1$ are
in set $M_{i1}$, and finally any lattice points off of $e_1$ are in
set $M_{i2}$. Hence the partition matrix $M_{ij}$ is exactly the
one constructed in the previous section. Since the polygons are
essential, at most one of them is $1$-dimensional, thus two of them,
say $P_0$ and $P_2$ as shown, have at least one point off the edge
$e_1$ marked ``2''.

\item[Case 2a]  The polygons are partitioned according to \rf{case2a}.

\begin{figure}[h]
\centerline{
\scalebox{0.60}{
\input{dim2resfig2.pstex_t}
             }
}           
  \caption{Case 2a}
         \label{F:case2a}
     \end{figure}

There are two distinguished edges $e_1$ and $e_2$.  Polygons
$P_0$ and $P_1$ share a complete flag along edge $e_1$ and are
partitioned accordingly into three sets $M_{i0}$, $M_{i1}$, and
$M_{i2}$ for $i = 0,1$ as shown.  But, this time the third polygon
$P_2$ has only one point on $e_1$ represented by the dotted line. This
point is put into set $M_{20}$ and all other points of $P_2$ are put
into set $M_{22}$. Notice that in the third polygon $M_{21} =\emptyset$.
Also note that $P_0$ and $P_1$ each have only one point (marked
``1'') on the dotted lines representing edges parallel to $e_2$.

To see that this is a partition matrix we must show that the sum of
lattice points in $M_{i0}$, $M_{j1}$, and $M_{k2}$ with $\{i,j,k\} =
\{0,1,2\}$ is in the interior of $P$. From the diagram this
corresponds to taking three points marked ``0'', ``1'' and ``2''
respectively from three different polygons and showing they cannot all
lie on parallel edges with the same inner normal.  If ``0'' and ``1'' come
from the first two polygons then their sum is in the interior of the
edge $e_1$ of the two-dimensional (by essentiality) sum of $P_0$ and
$P_1$. By definition however a point marked ``2'' from $P_2$ is not on
this edge.

If instead ``0'' comes from the $P_2$ and ``1'' comes from $P_0$ or $P_1$
the Minkowski sum of these two points is either the vertex lying only
on edges $e_1$ and $e_2$ or in the interior of $e_1$.  However, the
points marked ``2'' from the third polygon $P_0$ or $P_1$ are not on
either of these two edges.

For the combinatorial degree we apply \rl{singlecolor}.  This
shows that any face of $P$ (facet of $\hat{P}$) maximally canonically
colored by $C_2$ must only be colored by $C_2$. Such a face must have
its coloring matrix with $1$'s only in the first two columns. It is easy to
see from the picture that this can only happen for the bottom edge
$e_1$. It is also easy to see that the two vertices of this edge are
colored $\{1,2\}$ and $\{0,2\}$ respectively.  In particular
there is a unique complete flag colored $(\{2\}, \{1,2\})$, as
desired.

\item[Case 2b]  The difference between {Case 2a} and {Case 2b} is that $P_1$
now also contains an edge parallel to $e_2$. One attempt to
account for this would be to use the same partition as in {Case
2a} above except all the lattice points in $P_1$ along the edge $e_2$
are placed in $M_{11}$. 

If this were a partition matrix, edge $e_1$ would remain the only edge
colored just $\{2\}$ and its vertices would still be colored $\{1,2\}$
and $\{0,2\}$.

To check if this is a partition matrix, most of the arguments from the
previous case go through. If we take points marked ``0'' from $P_2$ and
``1'' from $P_0$ or $P_1$, the sum is either the vertex between edge
$e_1$ and $e_2$ or in the interior of one of these edges. Taking ``0''
from $P_1$ and ``1'' from $P_0$ again yields a point in the interior of
edge $e_1$.  However, if we take ``0'' from $P_0$ and ``1'' from $P_1$
there is a problem if there is some edge other than $e_1$ passing
through both these points.  But this can only happen if the edge $e_3$
of $P_0$ directly before $e_1$ passes through the endpoint of $e_2$ in
$P_1$, marked ``1'' as in \rf{case2b1}.

\begin{figure}[h]
\centerline{
\scalebox{0.60}{
\input{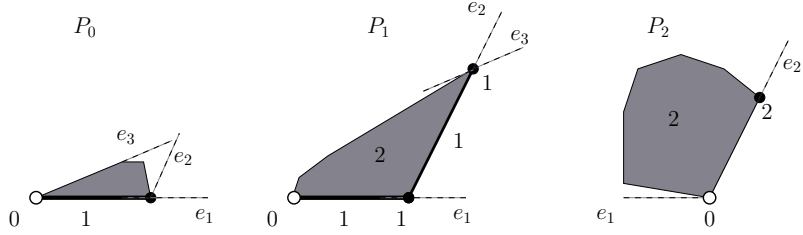}
             }
}           
  \caption{Case 2b: Failed partition}
         \label{F:case2b1}
     \end{figure}

In this case the partition fails, so we try to partition in a
different way, reversing the roles of $e_1$ and $e_2$ as in \rf{case2b2}.

\begin{figure}[h]
\centerline{
\scalebox{0.60}{
\input{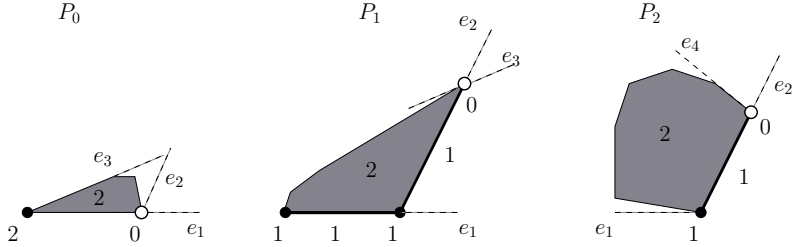}
             }
}           
  \caption{Case 2b: Switched partition}
         \label{F:case2b2}
     \end{figure}

Just as before this partition works as long as the next edge $e_4$ of
$P_2$ after $e_2$ does not pass through the left endpoint of $e_1$ in
$P_1$.  If both of the above attempts fail, we have that the edge
$e_3$ of $P_0$ before $e_1$ passes through the end point of $e_2$ in
$P_1$, and the edge $e_4$ of $P_2$ after $e_2$ passes through the left
endpoint of $e_1$ in $P_1$.  

In this final case, we show that we can find a partition matrix unless
we are in the exceptional case. First, assume that one of $P_0$
or $P_2$ is actually two-dimensional. Assume without loss of
generality it is $P_0$. Take the partition in \rf{case2b3}.

\begin{figure}[h]
\centerline{
\scalebox{0.60}{
\input{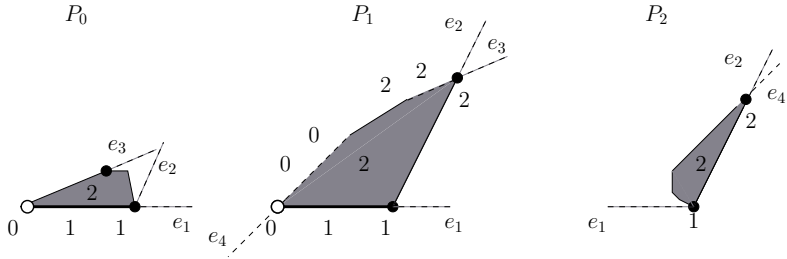}
             }
}           
  \caption{Case 2b:  Non-exceptional}
         \label{F:case2b3}
     \end{figure}

 Once again, a sum of points marked ``0'' from $P_0$ and ``1'' from $P_1$
must lie on the interior of $e_1$ which contains no points marked
``2'' from $P_2$. If we take ``0'' from $P_1$ then all the edges on which
it lies are between $e_1$ and $e_3$. But since $P_0$ was assumed
two-dimensional, the only such edge that could pass through a point
marked ``1'' is $e_1$ itself.

If we take ``1'' from $P_2$ we must lie on an edge of $P_2$
on or before $e_2$ but on or after $e_4$. The edge $e_2$ passes
only through points marked ``1'' or ``2'' in $P_0$ and $P_1$. The
next edge before $e_2$ is $e_1$ which passes through only ``0'' and
``1'' in $P_0$ and $P_1$. Furthermore all edges before $e_1$ and on
or after $e_4$ pass only through ``0'' in $P_1$ and, since $e_3$ is
before $e_4$, only through ``0'' in $P_0$. In every case we do not
get three points from different partition sets from three different
polygons all lying on the same edge.

A similar argument with an analogous partition applies if $P_2$ is
two-dimensional. Finally, if both $P_0$ and $P_2$ are one-dimensional,
the only way the partition above can fail is if the edge $e_3$,
parallel to $e_1$, which is known to go through a point marked ``2'' in
$P_1$, also goes through a point marked ``0''. But this can only happen
if this is the only other point of $P_1$ and moreover the edge
connecting this point to the endpoint of $e_1$ is parallel to $e_4$
which is also parallel to $e_2$. In other words, we must have $P_1$
have the same normal fan as the Minkowski sum of the two non-parallel
segments that are $P_0$ and $P_2$. That is to say we are in the
exceptional case.

\item[Case 3] The mixed case is somewhat easier and is partitioned as in \rf{mixed}.

\begin{figure}[h]
\centerline{
\scalebox{0.60}{
\input{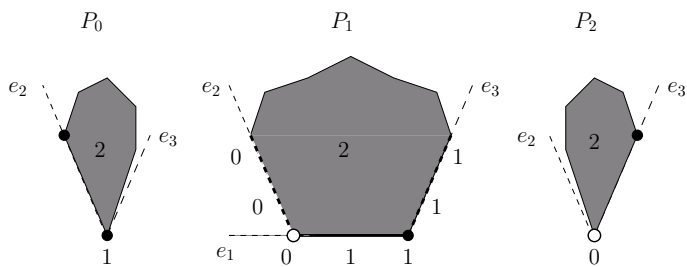}
             }
}           
  \caption{Case 3: Mixed partition}
         \label{F:mixed}
     \end{figure}

There are now three edges $e_1$, $e_2$, and $e_3$. There may actually
be several edges of $P_1$ between $e_2$ and $e_3$ in which case $e_1$
is the left most such edge. The edge $e_2$ intersects $P_1$ in either
a point or a whole edge, marked as a dashed line in the diagram. All
of the points on this edge are placed in set $M_{10}$. The edge $e_2$
together with $e_1$ and all other edges before $e_3$ intersects $P_2$
in a single point, represented by the dotted lines, placed in
partition set $M_{20}$ hence also marked ``0'' in the diagram. The
edge $e_3$ together with $e_1$ and all other edges before $e_3$
intersects $P_0$ in a single point placed in set $M_{01}$ Finally
$e_3$ intersects $P_1$ in either a single point or a whole edge marked
by a dashed line and all of the points are in set $M_{11}$. All the points
in $P_1$ on the edges between $e_2$ and $e_3$ are also in set $M_{11}$.
All other lattice points in each $P_i$ is in partition set $M_{i2}$ thus
marked as ``2''.

To show that this is a partition matrix we must take a point marked
``1'' from $P_0$ or $0$ from $P_2$ (otherwise we would have to take ``2''
from both). First if we took both ``0'' from $P_2$ and ``1'' from $P_0$,
then their Minkowski sum is on both $e_2$ and $e_3$. However any point
marked ``2'' from $P_1$ is on neither edge nor any edge in between. If
we took ``0'' from $P_2$ and ``1'' from $P_1$, the Minkowski sum is a
point lying only on edges on or between $e_1$ and $e_3$. However, any
point marked ``2'' from $P_0$ is not on any of these edges. The case
of ``1'' from $P_1$ is similar.

For combinatorial degree we note that the only edge colored just $\{2\}$
is the edge $e_1$.  Every other edge either intersects a point marked
``2'' or, if it is one of the other edges missing ``2'', can be colored
$\{0, 2\}$. The two vertices of this edge are colored $\{0,2\}$, and $\{1, 2\}$ respectively, completing the combinatorial
degree computation.
\end{description}
\end{pf}

To finish we show that an exceptional set of degrees {\em never}
has a compatible vertex partition resulting in an element
of residue one.

\begin{Prop} \label{P:exceptional} Let $\alpha_0, \alpha_1, \alpha_2$ be an exceptional
family of degrees on a two-dimensional toric variety $X$. Any
compatible vertex partition yields a residue matrix with 
determinant zero.
\end{Prop}

\begin{pf} We can assume that the corresponding polygons $P_0$ and
$P_1$ are one-dimensional and $P_2$ has the same normal fan as $P_0 +
P_1$.  Therefore $P_0$ and $P_1$ each have two vertices which we
denote by $u_0$, $u_1$ and $v_0$, $v_1$, and $P_2$ has four vertices $w_0$, $w_1$, $w_2$,
$w_3$.  For each pair $u_i$, $v_j$ there is a unique $w_k$ such that 
$u_i + v_j + w_k \in \inter(P_0 + P_1 + P_2)$. 
Now suppose we have a compatible vertex partition with associated partition matrix $M$.  
If $u_0$ and $u_1$ are in the same partition set, by the above we know
there do not exist $v_j$, $w_k$ such that both $u_0 + v_j + w_k$ and
$u_1 + v_j + w_k$ are in $\inter(P_0 + P_1 + P_2)$. Hence, there are
no non-zero terms in the expansion of the determinant of the induced
residue matrix. If $u_0$, $u_1$, $v_0$, $v_1$ all lie in two columns of the
partition matrix, then again there is no possible compatible choice of
$w_k$ in the complementary entry and the induced residue matrix will
have determinant zero. Therefore, up to relabeling we are left with only
one choice of partition matrix of the form:

$$\begin{bmatrix}
u_0 & u_1 & \emptyset \\
v_0 & \emptyset & v_1 \\
P_2^0 & P_2^1 & P_2^2 
\end{bmatrix}.
$$

For compatibility each of $P_2^i$ can only contain a unique
vertex $w_k$.  But this implies there is some vertex $w_k$ which
cannot lie in any of the $P_2^i$, a contradiction. Hence there are no
non-trivial compatible partition matrices. In particular there is no
partition yielding an element of residue one.
\end{pf}

This last result shows that while we know for essential degrees there
must always exist a polynomial of residue one by \rt{nonzero}, it cannot
always be obtained as the determinant of a matrix. On the flip side
the method does work for all but one quite degenerate situation.
We illustrate this by constructing residue matrices for some examples.

\begin{Example} Consider the polynomials:

\begin{align*}
f_0 &= a_0x + a_1xy + a_2y^2 \\
f_1 &= b_0 + b_1x + b_2x^2 + b_3xy \\
f_2 &= c_0 + c_1y +c_2xy^2.
\end{align*}

\begin{figure}[h]
\centerline{
\scalebox{0.60}{
\input{dim2example1.pstex_t}}}           
\caption{}
\label{F:dim2example1}
\end{figure}

The Newton polygons are shown in \rf{dim2example1} with the lattice
points labeled by their corresponding coefficients. This falls under
Case 3 of the previous theorem so applying the partition as in \rf{partition}
yields the following residue matrix:
$$\begin{bmatrix}
0 & a_0x & a_1xy + a_2y^2 \\
b_0 & b_1x +b_2x^2 & b_3xy \\
c_0 & 0 & c_1y + c_2xy^2.
\end{bmatrix}$$
The determinant is
$$ a_0b_3c_0xy - a_0b_0c_1xy  - a_0b_0c_2xy^2 - a_1b_1c_0x^2y - a_1b_2c_0x^3y - a_2b_1c_0xy^2 - a_2b_2c_0x^2y^2. $$
This is a polynomial supported on the interior of the Minkowski sum $P_0 + P_1 + P_2$. 
The homogenization up to critical degree has toric residue equal to 1.

\end{Example}

\begin{Example} Consider the polynomials:

\begin{align*} 
f_0 &= a_0 + a_1x \\
f_1 &= b_0 + b_1x + b_2y \\
f_2 &= c_0 + c_1xy. 
\end{align*}

\begin{figure}[h]
\centerline{
\scalebox{0.60}{
\input{dim2example2.pstex_t}}}           
\caption{}
\label{F:dim2example2}
\end{figure}

The corresponding Newton polygons are shown in \rf{dim2example2}.
These polygons can be classified under Case 2a of the above theorem.
As such we get the following residue matrix:
$$\begin{bmatrix}
a_0 & a_1x & 0 \\
b_0 & b_1x & b_2y \\
c_0 & 0 & c_1xy 
\end{bmatrix}.$$
The determinant is $$a_1b_2c_0xy + (a_0b_1c_1 - a_1b_0c_1)x^2y$$ 
which is supported in the interior of $P_0 + P_1 + P_2$ (consisting of two
points). Once again the homogenization up to critical degree yields
the desired element of residue 1.
\end{Example}

\begin{Example} \label{Ex:exbad}

Let us now consider an exceptional system:
\begin{align*}
f_0 &= a_0 + a_1x \\
f_1 &= b_0 + b_1x + b_2y + b_3xy\\
f_2 &= c_0 + c_1y.
\end{align*}
The Newton polygons consists of two line segments and their
Minkowski sum (a square). Since we are in the exceptional
case the theorem does not apply. 

However, there is a unique interior point of the Minkowski sum, so the
critical degree is trivial. Thus, there is a unique element of residue
$1$, namely the resultant itself which in this case is:
$$a_1b_0c_1 - a_0b_1c_1 - a_1b_2c_0 + b_3a_0c_0.$$
By \rp{exceptional} this polynomial is not expressible as the
determinant of a residue matrix.
\end{Example}

\section{Further Work and Conclusions} \label{S:conclusion}

Given a collection of $n+1$ semi-ample divisors on a toric variety $X$
which do not have a common zero, there exists a toric residue map
which is not identically zero if and only if the degrees of the
divisors are essential. The goal of this work was to explicitly construct
an element of residue one.

We have shown how compatible partitions of the Newton polytopes lead
to matrices whose determinant is an element of critical degree with
toric residue equal to a certain integer constant, namely the
combinatorial degree of a canonical induced coloring. In the case the
polytopes share a complete flag of faces and in almost every case in
dimension 2 we have shown how to choose this partition to yield an
element of residue exactly one.

The most obvious open question is to find compatible partitions
yielding elements of residue one in higher dimensions when the
polytopes do not necessarily share a complete flag. We have computed a
large number of examples in dimension three where the four polytopes
are simplices. In every case we have found working partitions. Of
course there will be exceptional families, as in dimension two, where
no such partitions exist. However, it is hoped that these will be
relatively rare and perhaps nonexistent in the most important case
when the polytopes are all full dimensional.


\section*{Acknowledgments}

This project was inspired by discussions with Eduardo Cattani and
Alicia Dickenstein. We would like to especially thank David Cox for
his careful proofreading of the text and all of his advice and
suggestions.


\end{document}